\documentclass[12pt]{amsart}
\usepackage{latexsym, amssymb, amsmath, amscd}
 \usepackage{eucal}

\usepackage[dvips]{graphicx}
\usepackage[latin1]{inputenc}
\usepackage[english]{babel}

    \newtheorem{rema}{Remark}[section]
    \newtheorem{propo}[rema]{Proposition}

   \newtheorem{theo}[rema]{Theorem}
   \newtheorem{def-theo}[rema]{Definition-Theorem}
 
   \newtheorem{defi}[rema]{Definition}
    \newtheorem{lemma}[rema]{Lemma}
    \newtheorem{corol}[rema]{Corollary}
     \newtheorem{exam}[rema]{Example}
  \newtheorem{rmk}[rema]{Remark}

	\newcommand{\nno}{\nonumber}

 \newcommand{\pf}{{\it Proof:}\hspace{2ex}}
 
 \newcommand{\epfv}{\hspace{1em}$\Box$\vspace{1em}}

\newcommand{\bQ}{{\mathbb Q}}
\newcommand{\bN}{{\mathbb N}}
\newcommand{\bT}{{\mathbb T}}
\newcommand{\bbT}{\bar{\mathbb T}}
\newcommand{\cT}{{\mathcal T}}
\newcommand{\cP}{{\mathcal P}}

\newcommand{\cS}{{\mathcal S}}

\newcommand{\cH}{{\mathcal H}}
\newcommand{\cV}{{\mathcal V}}

\newcommand{\cNcs}{{${\mathcal N}$CS} }
\newcommand{\cNsf}{{{\mathcal N}Sym}}
\newcommand{\cQf}{{{\mathcal Q}Sym}}

\newcommand{\klam}{{K\langle \Lambda \rangle}}

\newcommand{\BQ}{\begin{eqnarray}}
\newcommand{\EQ}{\end{eqnarray}}
\newcommand{\BQn}{\begin{eqnarray*}}
\newcommand{\EQn}{\end{eqnarray*}}
\newcommand{\BL}{\begin{align}}
\newcommand{\EL}{\end{align}}
\newcommand{\BLn}{\begin{align*}}
\newcommand{\ELn}{\end{align*}}
\newcommand{\BA}{\begin{align}}
\newcommand{\EA}{\end{align}}
\newcommand{\BAn}{\begin{align*}}
\newcommand{\EAn}{\end{align*}}
\newcommand{\wtilde}{\widetilde}

\newcommand{\lp}{\left(}
\newcommand{\rp}{\right)}

\newcommand{\tun}{\begin{picture}(5,0)(-2,-1)
\put(0,0){\circle*{2}}
\end{picture}}

\newcommand{\tdeux}{\begin{picture}(7,7)(0,-1)
\put(3,0){\circle*{2}}
\put(3,0){\line(0,1){5}}
\put(3,5){\circle*{2}}
\end{picture}}

\newcommand{\ttroisun}{\begin{picture}(15,8)(-5,-1)
\put(3,0){\circle*{2}}
\put(-0.65,0){$\vee$}
\put(6,7){\circle*{2}}
\put(0,7){\circle*{2}}
\end{picture}}
\newcommand{\ttroisdeux}{\begin{picture}(5,12)(-2,-1)
\put(0,0){\circle*{2}}
\put(0,0){\line(0,1){5}}
\put(0,5){\circle*{2}}
\put(0,5){\line(0,1){5}}
\put(0,10){\circle*{2}}
\end{picture}}

\newcommand{\tquatreun}{\begin{picture}(15,12)(-5,-1)
\put(3,0){\circle*{2}}
\put(-0.65,0){$\vee$}
\put(6,7){\circle*{2}}
\put(0,7){\circle*{2}}
\put(3,7){\circle*{2}}
\put(3,0){\line(0,1){7}}
\end{picture}}

\newcommand{\tquatrequatre}{\begin{picture}(15,18)(-5,-1)
\put(3,5){\circle*{2}}
\put(-0.65,5){$\vee$}
\put(6,12){\circle*{2}}
\put(0,12){\circle*{2}}
\put(3,0){\circle*{2}}
\put(3,0){\line(0,1){5}}
\end{picture}}
\newcommand{\tquatrecinq}{\begin{picture}(9,19)(-2,-1)
\put(0,0){\circle*{2}}
\put(0,0){\line(0,1){5}}
\put(0,5){\circle*{2}}
\put(0,5){\line(0,1){5}}
\put(0,10){\circle*{2}}
\put(0,10){\line(0,1){5}}
\put(0,15){\circle*{2}}
\end{picture}}

\newcommand{\tcinqun}{\begin{picture}(20,8)(-5,-1)
\put(3,0){\circle*{2}}
\put(-0.5,0){$\vee$}
\put(6,7){\circle*{2}}
\put(0,7){\circle*{2}}
\put(3,0){\line(2,1){10}}
\put(3,0){\line(-2,1){10}}
\put(-7,5){\circle*{2}}
\put(13,5){\circle*{2}}
\end{picture}}

\newcommand{\tcinqdix}{\begin{picture}(15,19)(-5,-1)
\put(3,5){\circle*{2}}
\put(-0.5,5){$\vee$}
\put(6,12){\circle*{2}}
\put(0,12){\circle*{2}}
\put(3,0){\circle*{2}}
\put(3,0){\line(0,1){12}}
\put(3,12){\circle*{2}}
\end{picture}}
\newcommand{\tcinqonze}{\begin{picture}(15,26)(-5,-1)
\put(3,5){\circle*{2}}
\put(-0.65,5){$\vee$}
\put(6,12){\circle*{2}}
\put(0,12){\circle*{2}}
\put(3,0){\circle*{2}}
\put(3,0){\line(0,1){5}}
\put(0,12){\line(0,1){7}}
\put(0,19){\circle*{2}}
\end{picture}}
\newcommand{\tcinqdouze}{\begin{picture}(15,26)(-5,-1)
\put(3,5){\circle*{2}}
\put(-0.65,5){$\vee$}
\put(6,12){\circle*{2}}
\put(0,12){\circle*{2}}
\put(3,0){\circle*{2}}
\put(3,0){\line(0,1){5}}
\put(6,12){\line(0,1){7}}
\put(6,19){\circle*{2}}
\end{picture}}
\newcommand{\tcinqtreize}{\begin{picture}(5,26)(-2,-1)
\put(0,0){\circle*{2}}
\put(0,0){\line(0,1){7}}
\put(0,7){\circle*{2}}
\put(0,7){\line(0,1){7}}
\put(0,14){\circle*{2}}
\put(-3.65,14){$\vee$}
\put(-3,21){\circle*{2}}
\put(3,21){\circle*{2}}
\end{picture}}
\newcommand{\tcinqquatorze}{\begin{picture}(9,26)(-5,-1)
\put(0,0){\circle*{2}}
\put(0,0){\line(0,1){5}}
\put(0,5){\circle*{2}}
\put(0,5){\line(0,1){5}}
\put(0,10){\circle*{2}}
\put(0,10){\line(0,1){5}}
\put(0,15){\circle*{2}}
\put(0,15){\line(0,1){5}}
\put(0,20){\circle*{2}}
\end{picture}}

\newcommand{\tddeux}[2]{\begin{picture}(12,5)(0,-1)
\put(3,0){\circle*{2}}
\put(3,0){\line(0,1){5}}
\put(3,5){\circle*{2}}
\put(6,-2){\tiny #1}
\put(6,3){\tiny #2}
\end{picture}}

\newcommand{\tdtroisun}[3]{\begin{picture}(20,12)(-5,-1)
\put(3,0){\circle*{2}}
\put(-0.65,0){$\vee$}
\put(6,7){\circle*{2}}
\put(0,7){\circle*{2}}
\put(5,-2){\tiny #1}
\put(9,5){\tiny #2}
\put(-5,5){\tiny #3}
\end{picture}}
\newcommand{\tdtroisdeux}[3]{\begin{picture}(12,12)(-2,-1)
\put(0,0){\circle*{2}}
\put(0,0){\line(0,1){5}}
\put(0,5){\circle*{2}}
\put(0,5){\line(0,1){5}}
\put(0,10){\circle*{2}}
\put(3,-2){\tiny #1}
\put(3,3){\tiny #2}
\put(3,9){\tiny #3}
\end{picture}}

\newcommand{\tdquatreun}[4]{\begin{picture}(20,12)(-5,-1)
\put(3,0){\circle*{2}}
\put(-0.6,0){$\vee$}
\put(6,7){\circle*{2}}
\put(0,7){\circle*{2}}
\put(3,7){\circle*{2}}
\put(3,0){\line(0,1){7}}
\put(5,-2){\tiny #1}
\put(8.5,5){\tiny #2}
\put(1,10){\tiny #3}
\put(-5,5){\tiny #4}
\end{picture}}

\newcommand{\tdquatretrois}[4]{\begin{picture}(20,20)(-5,-1)
\put(3,0){\circle*{2}}
\put(-.65,0){$\vee$}
\put(6,7){\circle*{2}}
\put(0,7){\circle*{2}}
\put(6,14){\circle*{2}}
\put(6,7){\line(0,1){7}}
\put(5,-2){\tiny #1}
\put(9,5){\tiny #2}
\put(-5,5){\tiny #4}
\put(9,12){\tiny #3}
\end{picture}}

\newcommand{\tdquatrequatre}[4]{\begin{picture}(20,14)(-5,-1)
\put(3,5){\circle*{2}}
\put(-.65,5){$\vee$}
\put(6,12){\circle*{2}}
\put(0,12){\circle*{2}}
\put(3,0){\circle*{2}}
\put(3,0){\line(0,1){5}}
\put(6,-3){\tiny #1}
\put(6,4){\tiny #2}
\put(9,12){\tiny #3}
\put(-5,12){\tiny #4}
\end{picture}}

\newcommand{\tdquatrecinq}[4]{\begin{picture}(12,19)(-2,-1)
\put(0,0){\circle*{2}}
\put(0,0){\line(0,1){5}}
\put(0,5){\circle*{2}}
\put(0,5){\line(0,1){5}}
\put(0,10){\circle*{2}}
\put(0,10){\line(0,1){5}}
\put(0,15){\circle*{2}}
\put(3,-2){\tiny #1}
\put(3,3){\tiny #2}
\put(3,9){\tiny #3}
\put(3,14){\tiny #4}
\end{picture}}

\title[A \cNcs  System over the Grossman-Larson Hopf Algebra]
{A Noncommutative Symmetric System over the Grossman-Larson Hopf Algebra
of Labeled Rooted Trees}

    \author{Wenhua Zhao}      
     \date{\today}
    
\begin{document}

\begin{abstract}
In this paper, we construct explicitly a noncommutative symmetric (${\mathcal N}$CS) system 
over the Grossman-Larson Hopf algebra of labeled rooted trees. By the 
universal property of the \cNcs system formed by the generating functions 
of certain noncommutative symmetric functions, we obtain a specialization 
of noncommutative symmetric functions by labeled rooted trees. 
Taking the graded duals, we also get a graded Hopf algebra homomorphism from 
the Connes-Kreimer Hopf algebra of labeled rooted forests to 
the Hopf algebra of quasi-symmetric functions. 
A connection of the coefficients of the third generating function 
of the constructed \cNcs system
with the order polynomials of rooted trees 
is also given and proved.  
\end{abstract}

   

 \bibliographystyle{alpha}
    \maketitle


\renewcommand{\theequation}{\thesection.\arabic{equation}}
\renewcommand{\therema}{\thesection.\arabic{rema}}
\setcounter{equation}{0}
\setcounter{rema}{0}
\setcounter{section}{0}

\section{\bf Introduction}
 \label{S1}

Let $K$ be any unital commutative 
$\bQ$-algebra and $A$ a 
unital associative but not necessarily 
commutative $K$-algebra. 
Let $t$ be a formal central parameter, 
i.e. it commutes with all elements of $A$, 
and $A[[t]]$ the $K$-algebra 
of formal power series 
in $t$ with coefficients in $A$.
A {\it \cNcs $($noncommutative symmetric$)$ system} 
over $A$ (see Definition \ref{Main-Def}) 
by definition is a $5$-tuple 
$\Omega\in A[[t]]^{\times 5}$ 
which satisfies the defining equations 
(see Eqs.\,$(\ref{UE-0})$--$(\ref{UE-4})$) 
of the NCSFs (noncommutative symmetric functions) 
first introduced and studied 
in the seminal paper \cite{G-T}. 
When the base algebra 
$K$ is clear in the context,
the ordered pair $(A, \Omega)$ 
is also called a {\it \cNcs system}.   
In some sense, a \cNcs  system over 
an associative $K$-algebra can be viewed 
as a system of analogs in $A$
of the NCSFs defined by 
Eqs.\,$(\ref{UE-0})$--$(\ref{UE-4})$.
For some general discussions on 
the \cNcs  systems, see \cite{GTS-I}.
For a family of \cNcs  systems  
over differential operator algebras 
and their applications to 
the inversion problem, 
see \cite{GTS-II} and \cite{GTS-III}.
For more studies on NCSFs, 
see \cite{T}, \cite{NCSF-II}, 
\cite{NCSF-III}, \cite{NCSF-IV}, 
\cite{NCSF-V} and \cite{NCSF-VI}.

One immediate but probably the most 
important example of the \cNcs systems 
is $(\cNsf, \Pi)$ formed 
by the generating functions of 
the NCSFs defined in \cite{G-T}
by Eqs.\,$(\ref{UE-0})$--$(\ref{UE-4})$ 
over the free $K$-algebra $\cNsf$ of NCSFs  
(see Section \ref{S2}). 
It serves as the universal \cNcs  system 
over all associative $K$-algebra 
(see Theorem \ref{Universal}).  
More precisely, for any \cNcs  system $(A, \Omega)$, 
there exists a unique $K$-algebra homomorphism 
$\cS: \cNsf \to A$ such that 
$\cS^{\times 5}(\Pi) = \Omega$ 
(here we have extended the homomorphism
$\cS$ to $\cS: \cNsf[[t]] \to A[[t]]$ 
by the base extension). 

The universal property of 
the \cNcs  system $(\cNsf, \Pi)$
can be applied as follows 
when a \cNcs  system $(A, \Omega)$ 
is given.
Note that, as an important topic in 
the symmetric function theory, 
the relations or polynomial identities 
among various NCSFs have been 
worked out explicitly (see \cite{G-T}). 
When we apply the $K$-algebra 
homomorphism $\cS:\cNsf \to A$ 
guaranteed by the universal property of 
the system $(\cNsf, \Pi)$ 
to these identities, 
they are transformed into 
identities among the corresponding
elements of $A$ in the system $\Omega$.
This will be a very effective way 
to obtain identities for 
certain elements of $A$ if 
we can show they are 
involved in a \cNcs  system over $A$.
On the other hand, 
if a \cNcs  system 
$(A, \Omega)$ has already been 
well-understood, 
the $K$-algebra homomorphism 
$\cS:\cNsf  \to A$ in return 
provides a {\it specialization} 
or {\it realization} (\cite{G-T}, \cite{St2}) 
of NCSFs, which may provide some 
new understandings on NCSFs.
For more studies on the specializations 
of NCSFs, see the references quoted 
above for NCSFs.

In this paper, we  apply 
the gadget above to the Grossman-Larson 
Hopf algebra of labeled rooted trees.
To be more precise, for any non-empty 
$W\subseteq \bN^+$,\footnote{All constructions and results 
of this paper work equally well for any non-empty weighted set $W$, 
i.e. any non-empty set $W$ with a fixed weight function 
$wt: W\to \bN^+$ such that, for any $k\in\bN^+$, $wt^{-1}(k)$ is a 
finite subset of $W$. But, for simplicity and convenience, 
we will always assume that $W$ is a non-empty subset of $\bN^+$.}
let $\cH_{GL}^W$ 
be the Grossman-Larson Hopf 
algebra $\cH_{GL}^W$ (\cite{GL}, \cite{CK}, \cite{F})
of rooted trees labeled 
by positive integers of $W$.
We first introduce five generating functions
(see Eqs.\,(\ref{Def-3f(t)})--(\ref{Def-3m(t)})
and Eq.\,(\ref{Def-3d(t)})) of certain elements 
of $\cH_{GL}^W$ and show that they form  
a \cNcs  system $\Omega_\bT^W$ 
over $\cH_{GL}^W$ 
(see Theorem \ref{Main-Thm-Trees}). 
Then, by the universal 
property of the \cNcs  system 
$(\cNsf, \Pi)$ from NCSFs,
we obtain a graded Hopf algebra 
homomorphism $\cT_W: \cNsf \to \cH_{GL}^W$, 
which gives a specialization  
of NCSFs by $W$-labeled 
rooted trees (see Theorem \ref{T4.6}).
By taking the graded duals, 
we get a graded Hopf algebra 
homomorphism $\cT_W^*: \cH_{CK}^W\to \cQf$ from 
the Connes-Kreimer Hopf algebra
$\cH_{CK}^W$ (\cite{Kr}), \cite{CK}, \cite{F}) to 
the Hopf algebra $\cQf$ of quasi-symmetric functions 
(\cite{Ge}, \cite{MR}, \cite{St2}).
Later in \cite{GTS-V}, 
it will be shown that, 
when $W=\bN^+$, the specialization 
$\cT_W$ above is actually an embedding and hence 
the Hopf algebra homomorphism $\cT_W^*$ is onto.
Finally, we give a
combinatorial interpretation of
the constants $\theta_T$ (see 
Definition \ref{Def-theta})
of rooted trees $T$ that appeared 
in the third component $\tilde d(t)$ 
(see Eq.\,$(\ref{Def-3d(t)})$) of the \cNcs  system 
$\Omega_\bT^W$ above.
We show that, for each rooted tree $T$,
the constant $\theta_T$ 
coincides with the coefficient of $s$ 
in the order polynomial $\Omega(B_-(T), s)$ 
(see \cite{St2}), where $B_-(T)$ is the rooted 
forest obtained by cutting off the root 
of $T$.

The arrangement of the paper is as follows.
In Section \ref{S2}, we first recall 
the definition of the \cNcs  systems (\cite{GTS-I}) 
$\Omega$ over the $K$-algebras $A$ 
and a result (see Proposition \ref{bialg-case}) 
on the \cNcs  systems when $A$  
is further a bialgebra or Hopf algebra.
We then recall the universal \cNcs  system 
$(\cNsf, \Pi)$ formed by generating functions 
of certain NCSFs in \cite{G-T}. 
The Hopf algebra structure 
of  $\cNsf$ and the universal property of 
the \cNcs  system $(\cNsf, \Pi)$ 
(see Theorem \ref{Universal})
will also be reviewed. 

In Section \ref{S3}, 
we first fix certain notation on 
rooted trees and recall 
the Connes-Kreimer Hopf algebra
$\cH_{CK}^W$ and 
the Grossman-Larson Hopf algebra 
$\cH_{GL}^W$ of $W$-labeled rooted forests 
and rooted trees, respectively.
Then, by using the duality between 
the Grossman-Larson Hopf algebra
and the Connes-Kreimer Hopf algebra 
(see Theorem \ref{Duality}), we prove a  
technic lemma, Lemma \ref{key-lemma-2}, 
that will be crucial for 
our later arguments. 

In Section \ref{S4}, 
we introduce five generating functions 
of $W$-labeled rooted trees 
and show in Theorem \ref{Main-Thm-Trees} 
that they form a \cNcs  system 
$\Omega_\bT^W$ over the Grossman-Larson 
Hopf algebra $\cH_{GL}^W$. 
By the universal property of 
the system $(\cNsf, \Pi)$, 
we get a graded $K$-Hopf algebra 
homomorphism ${\mathcal T}_W: \cNsf\to \cH_{GL}^W$
(see Theorem \ref{T4.6}).
By taking the graded duals, 
we get a graded Hopf algebra 
homomorphism $\cT_W^*: \cH_{CK}^W\to \cQf$ from 
the Connes-Kreimer Hopf algebra
$\cH_{CK}^W$ to the Hopf algebra $\cQf$ 
of quasi-symmetric functions 
(see Corollary \ref{cQsf}).

In Section \ref{S5}, we first recall 
the strict order polynomials
and the order polynomials of finite posets 
(partially ordered sets). 
Then, by applying some of results 
proved in \cite{WZ} for the strict 
order polynomials of rooted forests
and the well-known Reciprocity Relation 
(see Proposition \ref{Recip-Re}) between the 
strict order polynomials and the order polynomials
of finite posets, 
we show in Proposition \ref{S6-Main} 
that, for any $T\in \bar{\mathbb T}$,
the constant $\theta_T$ involved in the 
third component of the \cNcs  system 
$\Omega_\bT^W$ is same as the coefficient 
of $s$ of the order polynomial 
$\Omega(B_-(T), s)$ of the rooted forest 
$B_-(T)$. 

Finally, two remarks are as follows. First, 
as we pointed out early, by applying the specialization 
${\mathcal T}_W:\cNsf \to \cH_{GL}^W$,  
we will get a host of identities for the rooted trees involved 
in the \cNcs  system $\Omega_\bT^W$ from the identities 
of the NCSFs in $\Pi$. We believe some of these identities 
are interesting, at least from a combinatorial point view. 
But, to keep this paper in a certain size, we have to ask 
the reader who is interested to do the translations 
via the Hopf algebra homomorphism 
${\mathcal T}_W:\cNsf \to \cH_{GL}^W$. 
Secondly, some relations between 
the \cNcs  system $(\cH_{GL}^W, \Omega_\bT^W)$ 
constructed in this paper 
and the \cNcs  systems constructed in 
\cite{GTS-II} over differential 
operator algebras will be further 
studied in the followed paper \cite{GTS-V}. 
Some consequences of those relations 
to the inversion problem (\cite{BCW} and \cite{E}) 
and specializations of NCSFs 
will also be derived there. In particular, 
it will be shown that, with the label set 
$W=\bN^+$,  the $K$-Hopf algebra homomorphism 
${\mathcal T}_W:\cNsf \to \cH_{GL}^W$ 
in Theorem \ref{T4.6} is actually an embedding.

{\bf Acknowledgment}:
The author is very grateful to both referees
for many invaluable suggestions to improve the paper. 
In particular,  all the diagrams of trees in this paper 
are due to one of the referees
who very kindly sent the author 
all the latex commands.

\renewcommand{\theequation}{\thesection.\arabic{equation}}
\renewcommand{\therema}{\thesection.\arabic{rema}}
\setcounter{equation}{0}
\setcounter{rema}{0}

\section{\bf The Universal \cNcs  System from 
Noncommutative Symmetric Functions}
 \label{S2}

In this section, we first 
recall the definition 
of the \cNcs  systems (\cite{GTS-I})
over associative algebras 
and some of the NCSFs 
(noncommutative symmetric functions) 
first introduced and studied 
in the seminal paper \cite{G-T}. 
We then discuss the universal property 
of the \cNcs  system formed by the generating 
functions of these NCSFs. The main result that 
we will need later is Theorem \ref{Universal} which 
was proved in \cite{GTS-I}.
For some general discussions on 
the \cNcs  systems, see \cite{GTS-I}.
For more studies on NCSFs, 
see \cite{T}, \cite{NCSF-II}, 
\cite{NCSF-III}, \cite{NCSF-IV}, 
\cite{NCSF-V} and \cite{NCSF-VI}.

Let $K$ be any unital commutative 
$\bQ$-algebra\footnote{For the reader who is mainly interested in 
the combinatorial aspects of the main results 
of this paper, the base field 
$K$ throughout this paper can be safely chosen to be the field 
$\mathbb Q$ of rational numbers.}
 and 
$A$ any unital associative but not necessarily commutative 
$K$-algebra. Let $t$ be a formal central parameter, 
i.e. it commutes with all elements of $A$, and $A[[t]]$ 
the $K$-algebra of formal power series 
in $t$ with coefficients in $A$. 
First let us recall the following 
main notion of this paper.

\begin{defi} \label{Main-Def}  $(\text{\cite{GTS-I}})$ \, 
For any  unital associative $K$-algebra $A$, a $5$-tuple $\Omega=$ 
$( f(t)$, $g(t)$, $d\,(t)$, $h(t)$, $m(t) ) 
\in A[[t]]^{\times 5}$ is said 
to be a {\it \cNcs $($noncommutative symmetric$)$ system}
over $A$ 
if the following equations are satisfied.
\allowdisplaybreaks{
\begin{align}
&f(0)=1 \label{UE-0}\\
& f(-t)  g(t)=g(t)f (-t)=1, \label{UE-1}   \\
& e^{d\,(t)} = g(t), \label{UE-2} \\
& \frac {d g(t)} {d t}= g(t) h(t), \label{UE-3}\\ 
& \frac {d g(t)}{d t} =  m(t) g(t).\label{UE-4}
\end{align} }
\end{defi}

When the base algebra $K$ is clear in the context, we also call 
the ordered pair $(A, \Omega)$ a {\it \cNcs  system}. 
Since \cNcs  systems often come from generating functions 
of certain elements of $A$ that are under the consideration, 
the components of $\Omega$ will also be refereed as 
the {\it generating functions} of their coefficients. 

All $K$-algebras $A$ that we are going to work on in this paper
are $K$-Hopf algebras. We will freely use some standard notions 
and results from the theory of bialgebras and Hopf algebras, 
which can be found in the standard text books 
\cite{Abe}, \cite{Knu} and \cite{Mon}. For example, 
by {\it a sequence of divided powers} of a bialgebra 
or Hopf algebra $A$ we mean a sequence $\{c_n\,|\, n\ge 0\}$ 
of elements of $A$
such that, for any $n\ge 0$, we have 
\begin{align*}
\Delta c_n= \sum_{k\ge 0} c_{k} \otimes c_{n-k}.
\end{align*}

The following result proved in \cite{GTS-I} later
will be useful to us.

\begin{propo}\label{bialg-case}
Let $(A, \Omega)$ be a \cNcs  system as above. 
Suppose $A$ is further  a $K$-bialgebra. 
Then the following statements are equivalent. 
\begin{enumerate}
\item[$(a)$] The coefficients of $f(t)$ form a sequence of divided powers of $A$.

\item[$(b)$] The coefficients of $g(t)$ form a sequence of divided powers of $A$.

\item[$(c)$] One $(\text{hence also all})$ 
of \,$d(t)$, $h(t)$ and $m(t)$ has all 
its coefficients primitive in $A$.
\end{enumerate}
\end{propo}

In the following remark, we would like to point out 
a connection of the notion of \cNcs systems with 
the notion of combinatorial Hopf algebras, which was first introduced 
by M. Aguiar, N. Bergeron and F. Sottile 
in \cite{ABS}. 

\begin{rmk}\label{new-rmk}
First, as pointed out in Remark $2.17$ in \cite{GTS-I}, 
when $A$ is a graded 
and connected Hopf algebra, and one of 
the statements of 
Proposition \ref{bialg-case} holds, 
say statement $(b)$. Furthermore assume 
in this case that the coefficients 
of $t^m$ $(m\geq 0)$ of $g(t)$ 
are homogeneous with grading $m$. 
Then the data $(A, g(t))$ is equivalent 
to a combinatorial Hopf algebra structure 
on the graded dual Hopf algebra $A^*$ of $A$.
For more details of the equivalence above, 
see Remark $2.17$ in \cite{GTS-I}. 

Since all other components 
of $\Omega$ are completely determined by 
$g(t)$ $($see Lemma $2.5$ of \cite{GTS-I} or Theorem \ref{Universal} below$)$, 
the notion of \cNcs systems under the conditions above 
is also equivalent to the notion of combinatorial Hopf algebras. 
Therefore, from this point of view, the notion of \cNcs systems 
generalizes the notion of combinatorial 
Hopf algebras to associative $K$-algebras $A$, 
since, for \cNcs systems over $A$, 
$A$ does not have to be a 
bialgebra nor Hopf algebra and 
the equivalent conditions in 
Proposition \ref{bialg-case} do not 
have to hold either. 

On the other hand, 
we would also like to point out that 
the notion of \cNcs emphasizes 
the whole package of five 
generating functions of elements of $A$ 
instead of just one. In other words, 
it emphasizes solutions of the system of 
equations Eq.\,$(\ref{UE-0})$--$(\ref{UE-4})$.
Even though, once one of the components of $\Omega$, say $g(t)$ again, 
is fixed, the other four will be given by 
the corresponding universal polynomials 
of NCSFs in coefficients 
of $g(t)$ $($see Theorem \ref{Universal} below$)$, 
it is very often not trivial at 
all what values of these universal polynomials are, 
or in other words, it is still far away from clear 
how to write down directly and explicitly 
the other four components of $\Omega$.

The main aim of this paper is to construct 
explicitly a \cNcs system $\Omega_{\bT}^W$ 
over the Grossman-Larson Hopf algebra $\cH_{GL}^W$
of $W$-labeled trees without using 
any universal polynomials of NCSFs. 
Once we get the \cNcs system $\Omega_{\bT}^W$ 
explicitly, then, by Theorem \ref{Universal},
these universal polynomials 
of NCSFs will be transformed into identities 
of coefficients of the corresponding 
components of $\Omega_{\bT}^W$ 
$($see Remark \ref{R4.8}$)$. 
Another immediate consequence 
is that we also get a very 
``visualizable representation", 
or more formally, a specialization 
of NCSFs by $W$-labeled rooted trees, 
which in return could be useful for studying and 
understanding certain properties of NCSFs.
\end{rmk}

Next, let us recall some of the NCSFs 
first introduced and studied in (\cite{G-T}). 

Let $\Lambda=\{ \Lambda_m\,|\, m\geq 1\}$ 
be a sequence of noncommutative 
free variables and $\cNsf$ or $\klam$ 
the free associative algebra 
generated by $\Lambda$ over $K$.  
For convenience, we also set $\Lambda_0=1$.
We denote by
$\lambda (t)$ the generating function of 
$\Lambda_m$ $(m\geq 0)$, i.e. we set
\begin{align}
\lambda (t):= \sum_{m\geq 0} t^m \Lambda_m 
=1+\sum_{k\geq 1} t^m \Lambda_m.
\end{align}

In the theory of NCSFs (\cite{G-T}), 
$\Lambda_m$ $(m\geq 0)$ is 
the noncommutative analog 
of the $m^{th}$ classical (commutative) 
elementary symmetric function 
and is called the {\it $m^{th}$ 
$(\text{noncommutative})$ 
elementary symmetric function.}

To define some other NCSFs, we consider 
Eqs.\,$(\ref{UE-1})$--$(\ref{UE-4})$ 
over the free $K$-algebra $\cNsf$
with $f(t)=\lambda(t)$. The 
solutions for $g(t)$, $d\,(t)$, 
$h(t)$, $m(t)$ exist and are unique, 
whose coefficients will be the NCSFs 
that we are going to define.
Following the notation in \cite{G-T} 
and \cite{GTS-I}, we denote the resulted 
$5$-tuple by 
\begin{align}
\Pi:= (\lambda(t),\, \sigma(t),\, \Phi(t),\, \psi(t),\, \xi(t))
\end{align}
and write the last 
four generating functions of 
$\Pi$ explicitly as follows.

\allowdisplaybreaks{
\begin{align}
\sigma (t)&=\sum_{m\geq 0} t^m S_m,  \label{lambda(t)} \\
\Phi (t)&=\sum_{m\geq 1} t^m \frac{\Phi_m}m  \label{Phi(t)}\\
\psi (t)&=\sum_{m\geq 1} t^{m-1} \Psi_m, \label{psi(t)}\\
\xi (t)&=\sum_{m\geq 1} t^{m-1} \Xi_m.\label{xi(t)}
\end{align}}

Now, for any $m\geq 1$, 
we define $S_m$ to be the 
{\it $m^{th}$ $(\text{noncommutative})$ complete 
homogeneous symmetric function} and
$\Phi_m $ (resp.\,\,$\Psi_m$) 
the {\it $m^{th}$ power sum symmetric function 
of the second $($resp.\,\,first$)$ kind}. 
Note that, $\Xi_m$ $(m\geq 1)$  were denoted by $\Psi_m^*$ 
in \cite{G-T}. Due to Proposition \ref{omega-Lambda} 
below,  the NCSFs $\Xi_m$ $(m\geq 1)$
do not play an important role in the NCSF theory 
(see the comments in page $234$ in \cite{G-T}). 
But, in the context of some other problems, 
relations of $\Xi_m$'s 
with other NCSFs, especially, with $\Psi_m$'s, 
are also important (see \cite{GTS-III}, for example).
So here, following \cite{GTS-I}, 
we call $\Xi_m \in \cNsf$ $(m\geq 1)$ 
the {\it $m^{th}$ $(\text{noncommutative})$ 
power sum symmetric function of the third kind}.

The following two propositions proved in \cite{G-T} 
and \cite{NCSF-II} will be very useful 
for our later arguments.

\begin{propo}\label{bases}
For any unital commutative $\bQ$-algebra $K$,
the free algebra $\cNsf$ is freely generated
by any one of the families of the NCSFs 
defined above.
\end{propo}

\begin{propo}\label{omega-Lambda}
Let $\omega_\Lambda$ be the anti-involution of 
$\cNsf$ 
which fixes $\Lambda_m$ $(m\geq 1)$.
Then, for any $m\geq 1$, we have
\begin{align}
\omega_\Lambda (S_m)&=S_m,  \label{omega-Lambda-e1}\\
\omega_\Lambda (\Phi_m)&=\Phi_m,\label{omega-Lambda-e2} \\
\omega_\Lambda (\Psi_m)&=\Xi_m. \label{omega-Lambda-e3}
\end{align}
\end{propo}

Next, let us recall the following graded 
$K$-Hopf algebra structure 
of $\cNsf$. It has been shown in 
\cite{G-T} that $\cNsf$ is the universal enveloping algebra 
of the free Lie algebra generated 
by $\Psi_m$ $(m\geq 1)$. Hence, it has a $K$-Hopf  
algebra structure as all other universal enveloping algebras 
of Lie algebras do. Its co-unit $\epsilon:\cNsf \to K$,
 co-product $\Delta$ and 
 antipode $S$ are uniquely determined by 
\begin{align}
\epsilon (\Psi_m)&=0, \label{counit} \\
\Delta (\Psi_m) &=1\otimes \Psi_m +\Psi_m\otimes 1, \label{coprod}\\
S(\Psi_m) & =-\Psi_m,\label{antipode}
\end{align}
for any $m\geq 1$. 

Next, we introduce the {\it weight} of NCSFs 
by setting the weight of 
any monomial $\Lambda_{m_1}^{i_1} 
\Lambda_{m_2}^{i_2} \cdots \Lambda_{m_k}^{i_k}$
to be $\sum_{j=1}^k i_j m_j$. 
For any $m\geq 0$, we denote by $\cNsf_{[m]}$ 
the vector subspace of $\cNsf$ spanned 
by the monomials of $\Lambda$ 
of weight $m$. Then it is easy to see that 
\begin{align}\label{Grading-cNsf}
\cNsf=\bigoplus_{m\geq 0} \cNsf_{[m]}, 
\end{align}
which provides a grading for $\cNsf$. 

Note that, it has been shown in \cite{G-T}, 
for any $m\geq 1$, the NCSFs 
$S_m, \Phi_m, \Psi_m \in  \cNsf_{[m]}$. 
By Proposition \ref{omega-Lambda}, 
this is also true for the NCSFs $\Xi_m$'s.
By the facts above and 
Eqs.\,(\ref{counit})--(\ref{antipode}), 
it is also easy to check that, 
with the grading given in Eq.\,(\ref{Grading-cNsf}), 
$\cNsf$ forms a graded $K$-Hopf algebra. 
Its graded dual is given 
by the space $\cQf$ of quasi-symmetric functions, 
which were first introduced by I. Gessel \cite{Ge} 
(see \cite{MR} and \cite{St2} for more discussions).

Now we come back to our discussions on the \cNcs  systems. 
From the definitions of the NCSFs above, 
we see that $(\cNsf, \Pi)$ obviously forms a \cNcs  system.
More importantly, as shown in Theorem $2.1$ in \cite{GTS-I}, 
we have the following important theorem on 
the \cNcs  system $(\cNsf, \Pi)$. 

\begin{theo}\label{Universal}
Let $A$ be a $K$-algebra and $\Omega$ 
a \cNcs  system over $A$. Then, 

$(a)$ There exists a unique $K$-algebra homomorphism 
$\cS: \cNsf\to A$ such that 
$\cS^{\times 5} (\Pi)=\Omega$.

$(b)$ If $A$ is further a $K$-bialgebra $($resp.\,\,$K$-Hopf algebra$)$ 
and one of the equivalent statements in Proposition \ref{bialg-case} 
holds for the \cNcs  system $\Omega$, then $\cS: \cNsf\to A$ is also 
a homomorphism of $K$-bialgebras $($resp.\,\,$K$-Hopf algebras$)$.
\end{theo}

\begin{rmk}\label{Comm-Case}
By applying the similar arguments as in the 
proof of Theorem \ref{Universal}, 
or simply taking the quotient over 
the two-sided ideal generated by the commutators 
of $\Lambda_m$'s, it is easy to see that, 
over the category of commutative $K$-algebras, 
the universal \cNcs  system 
is given by the generating functions of
the corresponding classical 
$($commutative$)$ symmetric functions $($\cite{Mc}$)$.
\end{rmk}

\renewcommand{\theequation}{\thesection.\arabic{equation}}
\renewcommand{\therema}{\thesection.\arabic{rema}}
\setcounter{equation}{0}
\setcounter{rema}{0}

\section{\bf The Grossman-Larson Hopf Algebra and 
The Connes-Kreimer Hopf Algebra} \label{S3}

Let $K$ be any unital commutative $\bQ$-algebra
and $W$ a non-empty subset of positive 
integers. In this section, we first 
fix some notations 
for unlabeled rooted trees and 
$W$-labeled rooted trees
that will be used throughout this paper. 
We then recall the Connes-Kreimer Hopf algebra and
the Grossman-Larson Hopf algebra
of $W$-labeled forests and $W$-labeled 
rooted trees, respectively. 
Finally, by using the duality between 
the Grossman-Larson Hopf algebra
and the Connes-Kreimer Hopf algebra 
(see Theorem \ref{Duality}), we prove a  
technic lemma, Lemma \ref{key-lemma-2}, 
that will play an important role in 
our later arguments.

First, let us fix the following notation 
which will be used throughout the rest of this paper.

\vskip2mm

{\bf Notation:}

\vskip2mm

By a {\it rooted tree} we mean a finite
1-connected graph with one vertex designated as its {\it root.} 
For convenience, we also view the empty set $\emptyset$ 
as a rooted tree and call it the {\it emptyset} rooted tree.
The rooted tree with a single vertex 
is called the {\it singleton} 
and denoted by $\circ$.
There are natural ancestral relations between vertices.  
We say a vertex $w$ is a {\it child} of vertex $v$ 
if the two are connected by an
edge and $w$ lies further from the root than $v$.  
In the same situation, we say $v$ is the {\it parent} of $w$.  
A vertex is called a {\it leaf}\/ if it has no
children. 

Let $W\subseteq \bN^+$ be 
any non-empty subset of positive 
integers.
A {\it $W$-labeled rooted tree} 
is a rooted tree with each vertex labeled by 
an element of $W$. If an element $m \in W$ 
is assigned to a vertex $v$, 
then $m$ is called the {\it weight} 
of the vertex $v$. When we speak of isomorphisms 
between unlabeled (resp.\,\,$W$-labeled) rooted trees, 
we will always mean isomorphisms 
which also preserve 
the root (resp.\,\,the root and also the labels of vertices).
We will denote by $\mathbb T$ (resp.\,\,$\mathbb T^W$) 
the set of isomorphism classes 
of all unlabeled (resp.\,\,$W$-labeled) rooted trees. 
A disjoint union of any finitely many rooted trees 
(resp.\,\,$W$-labeled rooted trees) 
is called a {\it rooted forest} 
(resp.\,\,$W$-labeled {\it rooted forest}).
We denote by $\mathbb F$ (resp.\,\,$\mathbb F^W$) 
the set of unlabeled (resp.\,\,$W$-labeled) rooted forests.  


With these notions in mind, we establish the following notation.
\begin{enumerate}
\item For any rooted tree $T\in \bT^W$, we set the following notation:

\begin{itemize}
\item $\text{rt}_T$ denotes the root vertex of $T$ and $O(T)$
the set of all the children of $\text{rt}_T$. We set
$o(T)=|O(T)|$ (the cardinal number of the set $O(T)$). 

\item $E(T)$ denotes the set of edges of $T$.

\item $V(T)$ denotes the set of  vertices of $T$ and $v(T)=|V(T)|$.

\item $L(T)$ denotes the set of leaves of $T$ and $l(T)=|L(T)|$


\item For any $v\in V(T)$, we define the {\it height} of $v$ 
to be the number of edges in the (unique) geodesic connecting 
$v$ to $\text{rt}_{T}$.
The {\it height} of $T$ is defined to be the maximum of the heights of
its vertices.

\item For any $T\in \bT^W$ and $T\neq \emptyset$, 
$|T|$ denotes the sum of the weights of 
all vertices of $T$. 
When $T=\emptyset$, we set $|T|=0$.

\item For any $T\in \bT^W$, we denote 
by $\text{Aut}(T)$ the automorphism group 
of $T$ and $\alpha(T)$ the cardinal 
number of $\text{Aut}(T)$.

\end{itemize}

\item Any subset of $E(T)$ is called a {\it cut} of $T$. 
A cut $C\subseteq E(T)$ is said to be {\it admissible} 
if no two different edges of $C$ lie in 
the path connecting the root and a leaf. We denote by 
$\mathcal C(T)$ the set of all admissible cuts of $T$. 
Note that, the empty subset $\emptyset$  of $E(T)$
and $C=\{e\}$  
for any $e\in E(T)$ 
are always admissible cuts. We will identify any edge 
$e\in E(T)$ with the admissible cut $C:=\{e\}$ 
and simply say the edge $e$ itself is 
an admissible cut of $T$.

\item For any 
$T \in \bT^W$ with $T\neq \circ$, 
let $C\in \mathcal C(T)$ 
be an admissible cut of $T$ 
with $|C|=m\geq 1$.
Note that, after deleting 
the edges in $C$ from $T$, 
we get a disjoint union of $m+1$ 
rooted trees, 
say $T_0$, $T_1$, ..., $T_m$ 
with $\text{rt}(T)\in V(T_0)$.
We define $R_C(T)=T_0 \in \mathbb T^W$ 
and $P_C (T)\in \mathbb F^W$ 
the rooted forest formed by $T_1$, ..., $T_m$. 


\item For any disjoint admissible cuts 
$C_1$ and $C_2$, we say
``$C_1$ lies above $C_2$", and write $C_1\succ
C_2$, if $C_2 \subseteq E(R_{C_1}(T))$. 
This merely says that all edges of $C_2$ remain 
when we remove all edges of $C_1$ and 
$P_{C_1}(T)$.  Note that this relation
is not transitive. 
When we write $C_1\succ \cdots \succ
C_r\,$ for $C_1,\ldots, C_r\in \mathcal C(T)$,   
we will mean that $C_i \succ C_j$ 
whenever $i<j$. 

\item For any $T\in \bT^W$, we say $T$ 
is a {\it chain} if 
its underlying rooted tree
is a rooted tree with a single leaf.
We say $T$ is a {\it shrub} if 
its underlying rooted tree
is a rooted tree 
of height $1$. 
We say $T$ is {\it primitive} if 
its root has only one child. 
For any $m\geq 1$, we set $\mathbb H_m$, 
$\mathbb S_m$ and $\mathbb P_m$ to be the sets of 
the chains, shrubs and primitive rooted trees 
$T$ of weight $|T|=m$, respectively. 
$\mathbb H$, $\mathbb S$ and $\mathbb P$ 
are set to be the unions of $\mathbb H_m$, 
$\mathbb S_m$ and $\mathbb P_m$, 
respectively, for all $m\geq 1$.   

\end{enumerate}

For example, in the case where $W=\{1\}$, which allows not to write
the labels, we have 
\begin{eqnarray*}
\mathbb{H}&=&\left\{\tun,\tdeux,\ttroisdeux,\tquatrecinq,\tcinqquatorze\ldots \right\},\\
\mathbb{S}&=&\left\{\tdeux,\ttroisun,\tquatreun,\tcinqun\ldots \right\},\\
\mathbb{P}&=&\left\{\tdeux,\ttroisdeux,\tquatrequatre,\tquatrecinq,\tcinqdix,\tcinqonze,
\tcinqdouze,\tcinqtreize,\tcinqquatorze\ldots \right\}.
\end{eqnarray*}

For any non-empty $W\subseteq\bN^+$,
we define the following operations 
for $W$-labeled rooted forests.  For any 
$F\in \mathbb F^W$ which is disjoint 
union of $W$-labeled rooted trees 
$T_i$ $(1\leq i\leq m)$, 
we set $B_+(T_1, T_2, \cdots, T_m)$ \label{Tree-B+}
the rooted tree obtained by connecting roots 
of $T_i$ $(1\leq i\leq m)$ to a newly added root. 
We will keep the labels for the vertices of 
$B_+(T_1, T_2, \cdots, T_m)$ 
from $T_i$'s, but for the root, 
we label it by $0$. For convenience, 
we also fix the following short convention 
for the operation $B_+$. For the 
empty rooted tree $\emptyset$, we set
$B_+(\emptyset)$ to be the singleton 
labeled by $0$. For any $T_i\in \bT^W$  
$(1\leq i\leq m)$ and $j_i \geq 1$, 
the notation $B_+(T_1^{j_1}, T_2^{j_2}, 
\cdots, T_m^{j_m} )$
denotes the rooted tree obtained by applying 
the operation $B_+$ to $j_1$-copies of $T_1$; 
$j_2$-copies of $T_2$; and so on. Later,
for any unital $\bQ$-algebra $K$ and
$m\geq 1$, we will also
extend the operation $B_+$ 
multi-linearly to a linear map $B_+$ from
$\lp \cH_{CK}^W \rp^{\times m}$ to $\cH_{GL}^W$, 
where $\cH_{CK}^W$ and $\cH_{GL}^W$ at this moment 
are the vector spaces spanned over $K$ by 
the elements of $\bT^W$ and $B_+(\bT^W)$, 
respectively.

Now, we set $\bar {\mathbb T}^W\!:=\{ B_+(F) \, | \, F \in \mathbb F^W \}$.
Then,  $B_+\!: \mathbb F^W \to \bar{\mathbb T}^W$ becomes a bijection.
We denote by $B_- :  \bar{\mathbb T}^W \to \mathbb F^W$ 
the inverse map of $B_+$. More precisely, for any 
$T\in \bar{\mathbb T}^W$, 
$B_-(T)$ is the $W$-labeled rooted forest obtained by cutting off 
the root of $T$ as well as all edges connecting to the root in $T$.

Note that, precisely speaking,
elements of $\bar{\mathbb T}^W$ are not 
$W$-labeled trees for $0\not \in W$. 
But, if we set 
$\bar W=W\cup\{0\}$, 
then we can view $\bar{\mathbb T}^W$ 
as a subset 
of $\bar W$-labeled 
rooted trees $T$ with the root $\text{rt}_T$ 
labeled by $0$ 
and all other vertices 
labeled by non-zero elements of $\bar W$. 
We extend the definition of
the weight for elements 
of $\mathbb F^W$ to elements of 
$\bar{\mathbb T}^W$ by simply counting 
the weight of roots by zero. 
We set $\bar{\mathbb S}_m^W:=B_+(\mathbb S_m^W)$
$(m\geq 1)$ and $\bar{\mathbb S}^W:=B_+(\mathbb S^W)$. 
We also define
$\bar{\mathbb H}_m^W$, 
$\bar{\mathbb P}_m^W$, 
$\bar{\mathbb H}^W$ and 
$\bar{\mathbb P}^W$ 
in the similar way.

Next we fix a unital commutative $\bQ$-algebra $K$
and a non-empty subset of positive integers $W$, 
and first recall the Connes-Kreimer Hopf algebras 
$\mathcal H_{CK}^W$ of $W$-labeled rooted forests. 

As a $K$-algebra, the Connes-Kreimer Hopf algebra 
$\mathcal H_{CK}^W$ is the free commutative algebra 
generated by formal variables 
$\{X_T \,|\, T\in \mathbb T^W \}$. 
Here, for convenience, we will still use $T$ 
to denote the variable $X_T$ in 
$\mathcal H_{CK}^W$.
The $K$-algebra product is given by the disjoint union. 
The identity element of this algebra, denoted by $1$,  
is the free variable $X_\emptyset$ 
corresponding to the emptyset rooted tree $\emptyset$. 
The coproduct $\Delta: \mathcal H_{CK}^W \to 
\mathcal H_{CK}^W \otimes \mathcal H_{CK}^W$ 
is uniquely determined by setting 
\begin{align}
\Delta(1)&=1\otimes 1, \label{CK-Delta-1} \\
\Delta(T)&= T\otimes 1+ 
\sum_{C\in \mathcal C(T)} P_C(T) \otimes R_C(T). \label{CK-Delta-2}
\end{align}
The co-unit $\epsilon: \mathcal H_{CK}^W \to K$ 
is the $K$-algebra homomorphism 
which sends $1\in \mathcal H_{CK}^W$ to $1\in K$ 
and $T$ to $0$ for any $T\in \mathbb T^W$ 
with $T\neq \emptyset$. 
With the operations defined above and 
the grading given by the weight, 
the vector space $\mathcal H_{CK}^W$ 
forms a connected graded commutative 
bialgebra. Since any connected graded bialgebra 
is a Hopf algebra, there is a unique antipode 
$S: \mathcal H_{CK}^W\to \mathcal H_{CK}^W$ 
that makes $\mathcal H_{CK}^W$ a connected 
graded commutative $K$-Hopf algebra. 
For a formula for the antipode, see \cite{F}.

Next we recall the Grossman-Larson Hopf algebra 
of labeled rooted trees. As a vector space, 
the Grossman-Larson Hopf algebra $\mathcal H_{GL}^W$ 
is the vector space spanned by elements of $\bar{\mathbb T}^W$ 
over $K$. For any $T\in \bar{\mathbb T}^W$, we will still denote by 
$T$ the vector in $\mathcal H_{GL}^W$ that is corresponding to $T$. 
The algebra product is defined as follows. 

For any $T, S\in \bar{\mathbb T}^W$ with 
$T=B_+(T_1, T_2, \cdots, T_m)$, we set $T\cdot S$ to be the sum of  
the rooted trees obtained by connecting the roots of $T_i$ 
$(1\leq i\leq m)$ to vertices of $S$ 
in all possible $m^{v(S)}$ different ways. 
Note that, the identity element with respect to this 
algebra product is given by the singleton 
$\circ=B_+(\emptyset)$. But we will denote it by $1$.

To define the co-product $\Delta: \cH_{GL}^W \to \cH_{GL}^W \otimes \cH_{GL}^W$,
we first set 
\begin{align}
\Delta (\circ)=\circ \otimes \circ. 
\end{align}
 
Now let $T\in \bar{\mathbb T}^W$ with $T\neq \circ$, 
say $T=B_+ (T_1, T_2, \cdots, T_m)$ with
$m \geq 1$ and $T_i\in \mathbb T^W$ 
$(1\leq i\leq m)$. 
For any non-empty subset 
$I\subseteq \{1, 2, \cdots, m\}$,
we denote by 
$B_+(T_I)$ the rooted tree 
obtained by applying the $B_+$ operation 
to the rooted trees $T_i$ 
with $i\in I$. 
For convenience, when $I=\emptyset$, 
we set $B_+(T_I)=1$. 
With this notation fixed, the co-product 
for $T$ is given by
\begin{align}\label{GL-Delta}
\Delta (T)=\sum_{I\sqcup J=\{1, 2, \cdots, m\}} 
B_+(T_I) \otimes B_+(T_J).
\end{align}

The co-unit $\epsilon: \mathcal H_{GL}^W \to K$ 
is the $K$-algebra homomorphism 
which sends $1\in \mathcal H_{GL}^W$ to $1\in K$ 
and $T$ to $0$ for any $T\in \bbT^W$ 
with $T\neq \emptyset$. 
With the operations defined above and 
the grading given by the weight, 
the vector space $\mathcal H_{GL}^W$ 
forms a graded commutative 
bialgebra. Therefore, there is 
a unique antipode 
$S: \mathcal H_{CK}^W\to \mathcal H_{CK}^W$ 
that makes $\mathcal H_{CK}^W$ a 
graded $K$-Hopf algebra. Actually, 
by the general recurrent formula, 
we can write down the antipode 
of $\cH_{GL}^W$ as follows.

Note that the singleton $\circ$ is a group-like element 
and $S(\circ)=\circ^{-1}=\circ$. Now assume $T\neq \circ$ and write
$T=B_+(T_1, T_2, \cdots, T_m)$ with $m\geq 1$ and 
$T_i \in \bT^W$. Let $I:=\{1, 2,\cdots, m \}$. 
For $1\leq r\leq m$, let $\cP_r$ be the set of 
$r$-tuples $(I_1, I_2, \cdots, I_r)$ of disjoint 
non-empty subsets of $I$ whose union is $I$. 
In other words, $\cP_r$ is the set of all 
ordered partitions of $I$ 
into $r$ non-empty subsets 
of $I$.

\begin{lemma}
Let $S$ denote the antipode of the Grossman-Larson Hopf algebra $\cH_{GL}^W$
of $W$-labeled rooted trees. Then, for any $T\in \bbT^W$ with $T \neq \circ$ as above, 
we have 
\begin{align}\label{Anti-GL}
S(T)=\sum_{r=1}^m (-1)^r
\sum_{(I_1, \cdots, I_r) \in \cP_r }  B_+(T_{I_1})\, B_+(T_{I_2}) \cdots \, B_+(T_{I_r}) 
\end{align}
\end{lemma}
\pf By the general recurrent formula for the antipode of connected
cocommutative graded Hopf algebras, we know that the antipode $S$ 
of $\cH_{GL}^W$ satisfies the following equation:  
\begin{align}\label{Recurrent-GL}
S(T)=-T- \sum_{(I, J)\in \cP_2} 
S(B_+(T_I))\, B_+(T_J).
\end{align}
Then it is easy to check directly that, 
for any $T\in\bbT^W$, $S(T)$ given by Eq.\,(\ref{Anti-GL}) 
does satisfy Eq.\,(\ref{Recurrent-GL}). Since the solution to   
Eq.\,(\ref{Recurrent-GL}) is unique, the antipode $S$ of 
$\cH_{GL}^W$ is actually given by Eq.\,(\ref{Anti-GL}).
\epfv

Note that, from Eq.\,(\ref{GL-Delta}), 
it is easy to see that,
a rooted tree $T\in \bbT^W$ 
is a primitive element 
of the Hopf algebra $\cH_{GL}^W$ 
iff it is a primitive rooted tree 
in the sense that we defined before, 
namely the root of $T$ has 
one and only one child. It is noticeable that the set of 
primitive $W$-labeled rooted trees is a basis of the space  
$Prim(\cH_{GL}^W)$ of primitive 
elements of $\cH_{GL}^W$. Moreover, by the Milnor-Moore 
theorem, $\cH_{GL}^W$ 
is isomorphic to ${\mathcal U}(Prim(\cH_{GL}^W))$.

The relation between the Grossman-Larson Hopf algebra $\cH_{GL}^W$ and
the Connes-Kreimer Hopf algebra $\cH_{CK}^W$ is 
given by the following important theorem, 
which was proved in \cite{H} and \cite{F}.

\begin{theo}\label{Duality}
The Hopf algebras $\cH_{GL}^W$ 
and $\cH_{CK}^W$ are graded dual to each other. 
The pairing is given by, 
for any $T\in \bar{\mathbb T}^W$ 
and $S\in \mathbb F^W$, 
\begin{align}\label{pairing}
\langle T, F\rangle =\begin{cases} 0, & \text{ if } T \not \simeq B_+(F),\\
\alpha (T), &\text{ if } T \simeq B_+(F).
\end{cases}
\end{align}
\end{theo}

Furthermore, the following theorem on
the algebra structure constants 
of $\cH_{GL}^W$ was also proved in 
\cite{H} and \cite{F}.

\begin{theo}\label{L3.1.2}
For any $T', S\in \bar{\mathbb T}^W$, 
We have
\begin{align}\label{L3.1.2-e1}
T' \cdot S=\sum_{T\in \bar{\mathbb T}^W} 
\sum_{\substack {C\in \mathcal C(T)\\ B_+( P_C(T)) \sim T', \\  R_C(T)\sim S. }}
\frac{\alpha(T') \alpha(S)}{\alpha(T)} \,\,  T. 
\end{align}
\end{theo}

Note that, Eq.\,$(\ref{L3.1.2-e1})$ suggests that it is much more
convenient to work with the basis 
$\{\mathcal V_T:= T/ \alpha(T) \,| \, T\in \bar{\mathbb T}^W \}$
than the basis $\{ T \,| \, T\in \bar{\mathbb T}^W \}$. For example, 
in terms of $\mathcal V_T$, Eq.\,$(\ref{L3.1.2-e1})$ becomes 
\begin{align}\label{Rewriting}
\mathcal V_{T'} \cdot \mathcal V_S=\sum_{T\in \bar{\mathbb T}^W} 
\sum_{\substack {C\in \mathcal C(T)\\ B_+( P_C(T)) \sim T', \\ R_C(T)\sim S. }} \mathcal V_T. 
\end{align}

Finally, we extend Theorem \ref{L3.1.2} to a 
more general setting (see Lemma \ref{key-lemma-2} below). 
It can be viewed as a generalization 
of Lemma $2.8$ in \cite{WZ} which essentially is the case of 
Lemma \ref{key-lemma-2} when only primitive rooted trees 
are involved. First, let us fix the following notation.

Let $\vec{C}=(C_1, \ldots, C_r)\in \mathcal C(T)^{\times r}$ be a sequence 
of admissible cuts
with $C_1 \succ \cdots \succ C_r$.
We define a sequence of $T_{\vec{C},1}, \ldots , T_{\vec{C},r+1} \in \bar{\mathbb T}^W$
\label{T-Ci} 
as follows: we first set  $T_{\vec{C},1}= B_+ (P_{C_1}(T))$ and let
$S_1= R_{C_1}(T)$.  Note that 
$C_2, \ldots, C_r\in \mathcal C (S_1)$. 
We then set $T_{\vec C, 2}=B_+ (P_{C_2}(S_1))$ 
and  $S_2= R_{C_2}(S_1)$ and repeat this procedure 
until we get $S_{r}=R_{C_r}(S_{r-1})$ and then set 
$T_{\vec C, r+1}=S_r$.  
In the case that, each $C_i$ $(1\leq i\leq r)$
consists of a single edge, 
say $e_i\in E(T)$,  we simply denote $T_{\vec{C},i}$ 
by $T_{e_i}$.

Now we fix a positive integer $r$ and let
$y=\{ y_T^{(i)}\,|\, 1\leq i\leq r; \,, T\in \bar{\bT}^W\}$
be a collection of commutative formal variables.

\begin{lemma}\label{key-lemma-2}  
For any $r\geq 1$, we have, 
\allowdisplaybreaks{
\BQ\label{key-lemma2-e1}
&{}&
\sum_{\substack{(T_1,\ldots,T_r) \in (\bar{\bT}^W)^r }}
\left[ y_{T_1}^{(1)}\mathcal V_{T_1}\right]
\cdots \left[y_{T_r}^{(r)}\mathcal V_{T_{r}}\right] \\
&{}&\qquad\qquad \qquad =
\sum_{T\in \bar{\bT}^W }\,\,\,\,
\sum_{\substack{\vec{C}=(C_1,\ldots,C_{r-1}) \in {\mathcal C} (T)^{r-1}\\
C_1\succ \cdots \succ C_{r-1} }}y_{T_{\vec{C},1}}^{(1)}\cdots
y_{T_{\vec{C},r}}^{(r)}  \mathcal V_T. \nno 
\EQ}
\end{lemma}

\pf We denote the LHS of Eq.\,(\ref{key-lemma2-e1}) by $Q$ and write it as 
\begin{align*}
Q=\sum_{T \in \bar{\bT}^W} y_T \cV_T. 
\end{align*}
Then, by Eq.\,(\ref{pairing}) $y_T=\langle Q, B_-(T)\rangle$ for any $T\in \bar{\bT}^W$, 
where $B_-(T)$ is the forest obtained by deleting the root of $T$. 
So:
\begin{align*}
y_T&=\sum_{(T_1,\ldots, T_r)} y_{T_1}^{(1)} \cdots y_{T_r}^{(r)}
\langle \cV_{T_1} \ldots \cV_{T_r},B_-(T)\rangle \\
&= \sum_{(T_1,\ldots,T_r)} y_{T_1}^{(1)}\cdots y_{T_r}^{(r)}
\langle \cV_{T_1} \otimes (\cV_{T_2} \cdots  \cV_{T_r}), \Delta (B_-(T)) \rangle \\
&= \sum_{(T_1,\ldots,T_r)} y_{T_1}^{(1)}\cdots y_{T_r}^{(r)}
\langle \cV_{T_1} \otimes \cV_{T_2} \otimes 
(\cV_{T_3} \cdots \cV_{T_r}), (I\otimes \Delta) \circ \Delta(B_-(T)) \rangle, \\
\intertext{where $I$ is the identity map of $\cH_{CK}^W$. Repeating the process above:} 
&= \sum_{(T_1, \ldots,T_r)} y_{T_1}^{(1)}\cdots y_{T_r}^{(r)}
\langle \cV_{T_1} \otimes \cdots \otimes \cV_{T_r},(I^{\otimes (r-2)} 
\otimes \Delta) \circ \cdots \circ \Delta(B_-(T)) \rangle.
\end{align*}
One the other hand, by definition of the coproduct of ${\cH}_{CK}^W$ 
and definition of $\succ$, we have
\begin{align*}
& (I^{\otimes (r-2)} \otimes \Delta)\circ \cdots \circ (I\otimes \Delta) \circ \Delta(B_-(T)) \\
& \qquad \quad =\sum_{\substack{\vec{C}=(C_1,\ldots,C_{r-1}) \in {\mathcal C} (T)^{r-1}\\
C_1\succ \cdots \succ C_{r-1} }} B_-(T_{\vec{C},1})\otimes \ldots\otimes B_-(T_{\vec{C},r}).
\end{align*}
Therefore, we get
\begin{align*}
y_T=
\sum_{\substack{\vec{C}=(C_1,\ldots,C_{r-1}) \in {\mathcal C} (T)^{r-1}\\
C_1\succ \cdots \succ C_{r-1} }}y_{T_{\vec{C},1}}^{(1)}\cdots
y_{T_{\vec{C},r}}^{(r)}.
\end{align*}
\epfv

\renewcommand{\theequation}{\thesection.\arabic{equation}}
\renewcommand{\therema}{\thesection.\arabic{rema}}
\setcounter{equation}{0}
\setcounter{rema}{0}

\section{\bf A \cNcs  System over the Grossman-Larson Hopf Algebra 
$\mathcal H_{GL}^W$ of $W$-Labeled Rooted Trees}\label{S4}

In this section, for any non-empty $W\subseteq\bN^+$,
we construct a \cNcs  system $\Omega_{\bT}^W$ 
over the Grossman-Larson Hopf algebra $\cH_{GL}^W$.
First, let us introduce the following generating 
functions of certain elements of 
$\mathcal H^W_{GL}$, which will be 
the components of the \cNcs  system $\Omega_{\bT}^W$
corresponding to $f(t)$, $g(t)$, $h(t)$ and $m(t)$
according the notation in Definition \ref{Main-Def}.

\begin{align}
\tilde f(t):&=\sum_{T \in \bar{\mathbb S}^W} (-1)^{o(T)+|T|}  t^{|T|} {\mathcal  V_T}
=1 + \sum_{\substack{T \in \bar{\mathbb S}^W \\
T \neq \circ }} (-1)^{o(T)+|T|}  t^{|T|} {\mathcal  V_T}, 
\label{Def-3f(t)} \\ 
\tilde g(t):&= \sum_{T \in \bar{\mathbb T}^W} t^{|T|} {\mathcal  V_T}
=1 + \sum_{\substack{ T  \in \bar{\mathbb T}^W \\
T \neq \circ }} t^{|T|} {\mathcal  V_T}, 
\label{Def-3g(t)}\\ 
\tilde h(t):&= \sum_{T \in \bar{\mathbb H}^W} t^{|T|-1}  \beta_T \mathcal  V_T, 
\label{Def-3h(t)} \\ 
\wtilde m(t):&= \sum_{T \in \bar{\mathbb P}^W}  t^{|T|-1} \gamma_T {\mathcal  V_T}, 
\label{Def-3m(t)}
\end{align} 
where,  for any $T\in \bar{\mathbb H}^W$ (resp.\,\,$T\in \bar{\mathbb P}^W$),
$\beta_T$ (resp.\,\,$\gamma_T$) is the weight 
of the unique leaf 
(resp.\,\,the unique child of the root) 
of $T$. Note that, for the singleton $T=\circ$, we have
$\beta_T=\gamma_T=0$. So $\tilde h(0)=\wtilde m(0)=0$.

For example, when $W=\{1, 2 \}$, we have
\begin{align*}
& \tilde{f}(t)=1+\cV_{\tddeux{}{$1$}}t+\left( \cV_{\tdtroisun{}{$1$}{$1$}}- \cV_{\tddeux{}{$2$}}\right) t^2
+\left( \cV_{\tdquatreun{}{$1$}{$1$}{$1$}}- \cV_{\tdtroisun{}{$2$}{$1$}} \right) t^3+ \cdots\\
& \tilde{g}(t)=1 + \cV_{\tddeux{}{$1$}}t + \left(  \cV_{\tddeux{}{$2$}} + \cV_{\tdtroisun{}{$1$}{$1$}} 
  + \cV_{\tdtroisdeux{}{$1$}{$1$}}\right)  t^2 \\ 
& \quad +\left(  \cV_{\tdtroisun{}{$2$}{$1$}} 
 + \cV_{\tdtroisdeux{}{$1$}{$2$}}+\cV_{\tdtroisdeux{}{$2$}{$1$}}+
 \cV_{\tdquatreun{}{$1$}{$1$}{$1$}}+\cV_{\tdquatretrois{}{$1$}{$1$}{$1$}}
+\cV_{\tdquatrequatre{}{$1$}{$1$}{$1$}}+\cV_{\tdquatrecinq{}{$1$}{$1$}{$1$}} \right)t^3+\cdots\\
& \tilde{h}(t)=\cV_{\tddeux{}{$1$}}+
\left( \cV_{\tdtroisdeux{}{$1$}{$1$}}+ 2\cV_{\tddeux{}{$2$}} \right) t
+\left( \cV_{\tdquatrecinq{}{$1$}{$1$}{$1$}}+ 
2\cV_{\tdtroisdeux{}{$1$}{$2$}}+\cV_{\tdtroisdeux{}{$2$}{$1$}}       \right)t^2+ \cdots \\
& \wtilde{m}(t)=\cV_{\tddeux{}{$1$}}+\left( \cV_{\tdtroisdeux{}{$1$}{$1$}}+ 2\cV_{\tddeux{}{$2$}} \right) t
+\left(\cV_{\tdtroisdeux{}{$1$}{$2$}}+2\cV_{\tdtroisdeux{}{$2$}{$1$}}   + \cV_{\tdquatrequatre{}{$1$}{$1$}{$1$}}+\cV_{\tdquatrecinq{}{$1$}{$1$}{$1$}}\right)t^2+\cdots\
\end{align*}

Note that, from Eq.(\ref{Def-3f(t)}), we have 
\begin{align}
\tilde f(-t)=\sum_{T \in \bar{\mathbb S}^W} (-1)^{o(T)}  t^{|T|} {\mathcal  V_T}
=1 + \sum_{\substack{T \in \bar{\mathbb S}^W \\
T \neq \circ }} (-1)^{o(T)}  t^{|T|} {\mathcal  V_T}. 
\label{Def-3f(-t)} 
\end{align}

A different way to look at the generating function 
$\tilde f(-t)$ is as follows.

For any $m\in W$, let $\kappa_m$ denote the singleton labeled by $m$ and set
\begin{align}\label{Def-kappa}
\kappa(t):=\sum_{m\in W} t^m \kappa_m.
\end{align}

\begin{lemma} \label{Taylor-for-3f(t)}
\begin{align}\label{Taylor-for-3f(t)-e1}
\tilde f(-t)=1+\sum_{d\geq 1} \frac{(-1)^{d}}{d!} B_+( \, \kappa (t)^d \,),
\end{align}
where $B_+( \, \kappa (t)^d \,)$ denotes the element obtained 
by applying $B_+$ to $d$-copies of $\kappa (t)$.
\end{lemma}

\pf First, it is easy to see that, the
only terms that can appear in the expansion of 
the RHS of Eq.\,(\ref{Taylor-for-3f(t)-e1}) are shrubs. 
Secondly, from the definition of the operation $B_+$, 
we see that $B_+$ is symmetric and multi-linear
in its components. Therefore, 
we can expand the term 
$B_+( \, \kappa (t)^d \,)$ into a linear 
combination of rooted trees $S\in \bbT^W$ 
in a similar way as we expand 
the power $( \sum_{m\in W} t^m u_m)^d$ 
for some free
commutative variables $u_m$ $(m\in W)$.

Now, for any shrub $S\in \bar{\mathbb S}^W$ 
with $S\neq \circ$, let $\{m_j\in W\,|\, 1\leq j\leq N \}$
be the set of all labels of the leaves of $S$.
Let $i_j\geq 1$ $(1\leq j\leq N)$ 
be the number of the $m_j$-labeled 
leaves of $S$. Then we have
\begin{align}
o(S)&=\sum_{1\leq j\leq N} i_j , \label{Taylor-for-3f(t)-pe1} \\
|S|&= \sum_{1\leq j\leq N} i_j m_j, \label{Taylor-for-3f(t)-pe2}\\
\alpha(S)&=\prod_{1\leq j\leq N} (i_j)!.\label{Taylor-for-3f(t)-pe3}
\end{align}

Now let us consider the coefficient $c_S$ of
$S$ in the linear expansion of the RHS of 
Eq.\,(\ref{Taylor-for-3f(t)-e1}). 
By the observations in the 
first paragraph of the proof and 
Eqs.\,$(\ref{Taylor-for-3f(t)-pe1})$-$(\ref{Taylor-for-3f(t)-pe3})$, 
it is easy to see that we have
\begin{align*}
c_S &= \frac {(-1)^{o(S)}t^{|S|}}{o(S)!} 
\binom {o(S)} {i_1, \cdots, i_N}  \\
&=(-1)^{o(S)} \frac {t^{|S|}}{\prod_{1\leq j\leq N} (i_j)! }\\
&=(-1)^{o(S)}\frac {t^{|S|}}{\alpha (S)},
\end{align*}
which is same as the coefficient of $S$ in $\tilde f(-t)$ 
since $\mathcal V_S=\frac 1{\alpha(S)} S$. Hence we are done.
\epfv

To define the generating function $\tilde d(t)$ 
for the third component of the under-construction 
\cNcs  system $\Omega_{\bT}^W$, we first need 
the following definition. 

\begin{defi}\label{Def-theta}
$(a)$ We define a constant $\theta_T\in \bQ$ 
for each unlabeled rooted tree $T$ as follows.

\begin{enumerate}
\item For the singleton $\circ$ and 
any non-primitive rooted tree $T\in \mathbb T$, 
i.e. $o(T)>1$, we set $\theta_\circ = \theta_T=0$.

\item  For $T=B_+(\circ)$, we set $\theta_T=1$.

\item  For any primitive $T\in \mathbb P$ with $v(T)\geq 3$, we define $\theta_T$ 
inductively by 
\begin{align}\label{theta-Recur}
\theta_T = 1- \sum_{m \geq 2} \frac 1{m!} 
\sum_{\substack{\vec e=(e_1,\ldots,e_{m-1})\in E(T)^{m-1}\\e_1\succ
\cdots\succ e_{m-1}}}
\theta_{T_{e_1}} \theta_{T_{e_2}} \cdots  \theta_{T_{e_{m}}}, 
\end{align}
\end{enumerate}
where $T_{e_i}$'s in the equation above 
have been defined before Lemma \ref{key-lemma-2}.

$(b)$ For any $W\subseteq \bN$ and $W$-labeled 
rooted tree $T$, we set $\theta_T:=\theta_{\bar T}$, 
where  ${\bar T}$
is the underlying unlabeled rooted tree of $T$.
\end{defi}

\begin{rmk}\label{rmk-theta}
As we will show later in Section \ref{S5},
the constant $\theta_T$ $(T\in \mathbb T)$ 
has a natural combinatorial interpretations as follows.
If we write $T=B_+(F)$ for 
some rooted forest $F\in \mathbb F$ and let $\Omega(F, s)$ 
be the order polynomial $(\text{see \,\cite{St1}})$
of $F$, then $\theta_T$ will be 
the coefficient of $s$ of $\Omega(F, s)$ 
$($see Proposition \ref{S6-Main} in Section \ref{S5}$)$.
Furthermore, if we denote by 
$\nabla: K[s]\to K[s]$ the linear operator 
which maps any $f(s)\in K[s]$ to $f(s)-f(s-1)$,  
then $\theta_T$ is also the coefficient of $s$ 
of the polynomial $\nabla \Omega(T, s)$ 
$($see Corollary \ref{C6.10}$)$.
In order to keep our on-going 
arguments more focus, 
we will postpone to Section \ref{S5} a detailed 
discussion on these combinatorial 
interpretations of $\theta_T$ 
$(T\in \mathbb T)$. 
\end{rmk}

\begin{exam} 
By Eqs.\,$(\ref{Exam6.2-e2})$, $(\ref{Exam6.2-e4})$, $(\ref{Def-varphi-P})$ 
and $(\ref{S6-Main-e1})$ in Section \ref{S5}, it is easy to check that, 
for the chains $C_m$'s and $B_+(S_m)$ of the shrubs $S_m$'s, we have
\begin{align}
\theta_{C_m}&=\frac 1{m-1} \quad \text{ for any $m\geq 2$. } \\
\theta_{B_+(S_m)}&=(-1)^{m} b_m \quad \text{ for any $m\geq 0$, } 
\end{align}
where $b_m$ $(m\geq 0)$ are the
Bernoulli numbers which are defined by the generating function
\begin{align}
\frac{x}{e^{x}-1}=\sum_{m=0}^{\infty}b_m \frac{x^{m}}{m!}.
\end{align}
\end{exam}

Now, we introduce the following generating function:
\begin{align}
\tilde d(t):= \sum_{T \in \bar{\mathbb P}^W}  t^{|T|} \theta_T {\mathcal  V_T}. 
\label{Def-3d(t)}
\end{align}

For example, when $W=\{1, 2 \}$, we have
\begin{align*}
\tilde{d}(t)=&\cV_{\tddeux{}{$1$}}t +\left( \frac 12 \cV_{\tdtroisdeux{}{$1$}{$1$}}+ 
\cV_{\tddeux{}{$2$}} \right) t^2 \\
&\quad +\left(\frac 12 \cV_{\tdtroisdeux{}{$1$}{$2$}}+\frac 12 \cV_{\tdtroisdeux{}{$2$}{$1$}}   + \frac 16 \cV_{\tdquatrequatre{}{$1$}{$1$}{$1$}}+ \frac 13\cV_{\tdquatrecinq{}{$1$}{$1$}{$1$}}\right)t^3+\cdots
\end{align*}

Set 
\begin{align}\label{Def-3-Omega}
\Omega_\bT^W:=
(\, \tilde f(t),\, \tilde g(t), 
\, \tilde d\,(t), \, \tilde h(t), \wtilde m(t)\, ). 
\end{align}
Then, the main result of this section is the following theorem.

\begin{theo}\label{Main-Thm-Trees}
For any non-empty set $W\subseteq\bN$, 
$\Omega_\bT^W$ forms a \cNcs  system over 
the Grossman-Larson Hopf algebra $\mathcal H^W_{GL}$. 
\end{theo}

\pf Note that, by Eqs.\,$(\ref{Def-3f(t)})$ and $(\ref{Def-3g(t)})$, 
we have $\tilde f(0)=\tilde g(0)=1$, hence it will be enough to show
Eqs.\,(\ref{UE-1})--(\ref{UE-4}) in Definition \ref{Main-Def}
are satisfied by the generating functions in $\Omega_\bT^W$.

Let us start with Eq.\,(\ref{UE-1}). First, note that, 
since $\tilde g(0)=1$,  $\tilde g(t)$ 
as an element of $\mathcal H_{GL}^W[[t]]$ 
does have both left and right inverses. 
So we only need show $\tilde f(-t) \tilde g(t)=1$ for 
$\tilde g(t)\tilde f(-t)=1$ will follow automatically. 
Secondly, by Eq.\,$(\ref{Def-3f(-t)})$ and Lemma \ref{key-lemma-2} 
with $y^{(1)}_T=(-1)^{o(T)} t^{|T|}$ if $T\in \bar{\mathbb S}^W$ and 
$0$ otherwise, and $y^{(2)}_T= t^{|T|}$ 
for any $T\in \bbT^W$, we have
\allowdisplaybreaks{
\begin{align}\label{Main-Thm-Trees-pe1}
\tilde f(-t)\tilde g(t) &= \lp  \sum_{T' \in \bar{\mathbb S}^W} (-1)^{o(T')} t^{|T'|} 
\mathcal  V_{T'} \rp \lp  \sum_{T'' \in \bar{\mathbb T}^W} t^{|T''|} \mathcal  V_{T''} \rp \\
&=\sum_{T \in \bar{\mathbb T}^W} t^{|T|} 
\lp \sum_{\substack{C \in \mathcal C(T)\\ B_+(P_C(T)) \in \bar{\mathbb S}^W}} 
(-1)^{o(B_+(P_C(T)))} \rp  \mathcal  V_T \nno
\end{align}}

First, note that, for any rooted tree $T$ and an 
admissible cut $C$ of $T$, $B_+(P_C(T))$
is a shrub iff each edge in $C$ is the unique edge 
connecting with a leaf of $T$.
Therefore, the set of all admissible cuts 
$C$ such that $B_+(P_C(T)) \in \bar {\mathbb S}^W$ 
is in $1$-$1$ correspondence with the set of subsets of leaves of $T$.
Secondly, when $B_+(P_C(T)) \in \bar {\mathbb S}^W$ 
for an admissible cut $C$, $o (B_+(P_C(T)) )$ is same 
as the cardinal number $|C|$ of the cut $C$. 
With these observations, for any $T\in \bbT^W$ 
with  $l(T):=|L(T)|>0$, we have
\begin{align}\label{Main-Thm-Trees-pe2}
\sum_{\substack{C \in \mathcal C(T)\\ B_+(P_C(T))  \in \bar{\mathbb S}^W}} (-1)^{|C|}
=\sum_{k=0}^{l(T)} (-1)^k \binom {l(T)}k=0.
\end{align}
Combining Eqs.\,$(\ref{Main-Thm-Trees-pe1})$ and $(\ref{Main-Thm-Trees-pe2})$, 
we get $\tilde f(-t)\tilde g(t)=1$ and hence Eq.\,(\ref{UE-1}) for 
the system $\Omega_\bT^W$.

Now, let us prove Eq.\,(\ref{UE-2}) as follows.
\allowdisplaybreaks{
\begin{align*}
e^{\tilde d\,(t)}&=\sum_{k\geq 0} \frac 1{k!} \tilde d\,(t)^k\\
&=1+ \sum_{k\geq 1} \frac 1{k!} 
(\sum_{T\in \bar{\mathbb P}^W}  t^{|T|} \theta_T \mathcal  V_T)^k\\
\intertext{Applying Lemma \ref{key-lemma-2}:}
&=1+ \sum_{T\in \bar{\mathbb P}^W} t^{|T|} \theta_T \mathcal  V_T \\
&\quad \quad + 
\sum_{k\geq 2} \frac 1{k!} \sum_{T \in \bbT^W} t^{|T|} \lp
\sum_{\substack{\vec{e}=(e_1,\ldots,e_{k-1} )\in E(T)^{k-1} \\e_1\succ \cdots
\succ e_{k-1} }}  \theta_{T_{\vec{e},1}} \cdots
\theta_{T_{\vec{e}, k}}  \rp \mathcal  V_T   \\
&=1+ \sum_{T \in \bar{\mathbb T}^W} t^{|T|} 
\lp  \theta_T +  \sum_{k\geq 2} \frac 1{k!} 
\sum_{\substack{\vec{e}=(e_1,\ldots,e_{k-1} )\in E(T)^{k-1} \\e_1\succ \cdots
\succ e_{k-1} }}\theta_{T_{\vec{e},1}} \cdots
\theta_{T_{\vec{e}, k}} \rp   \mathcal  V_T \\
\intertext{Applying Eq.\,(\ref{theta-Recur}):}
&=1+ \sum_{T\in \bar{\mathbb T}^W} t^{|T|} \mathcal  V_T\\
&=\tilde g(t).
\end{align*} }

Therefore, we get $e^{\tilde d\,(t)}=\tilde g(t)$, which is Eq.\,(\ref{UE-2})
for the system $\Omega_\bT^W$.

To prove Eq.\,(\ref{UE-3}), first, by Lemma \ref{key-lemma-2}, we have 
\begin{align}
\wtilde m(t) \tilde g(t) &= \lp  \sum_{T' \in \bar{\mathbb P}^W} \gamma (T') t^{|T'|-1} 
\mathcal  V_{T'} \rp \lp  \sum_{T'' \in \bar{\mathbb T}^W} t^{|T''|} \mathcal  V_{T''} \rp \\
&=\sum_{T \in \bar{\mathbb T}^W} t^{|T|-1} 
\lp \sum_{\substack{C \in \mathcal C(T)\\ B_+(P_C(T)) \in \bar{\mathbb P}^W}}
 \gamma (B_+(P_C(T)))\rp \mathcal  V_T \nno \\
&=\sum_{T \in \bar{\mathbb T}^W} t^{|T|-1} 
\lp \sum_{e \in E(T)} \gamma (B_+(P_e(T))) \rp
\mathcal  V_T, \nno
\end{align}
where the last equality follows from the fact that, 
for any $C \in \mathcal C(T)$, $o\lp B_+(P_C(T))\rp =1$ iff 
$C$ consists of a single edge.

Note that, for any $e \in E(T)$,  $\gamma (B_+(P_e(T)) )=wt(v'_e)$, 
where $v'_e$ is the vertex of $e$ which is further away 
from the root of $T$. Therefore, continuing with the equation above, 
we have
\begin{align}
\wtilde m(t) \tilde g(t) &=\sum_{T \in \bar{\mathbb T}^W}  |T| \, t^{|T|-1}  \mathcal  V_T
=\frac {d \tilde g(t)}{dt}.\nno
\end{align}
Hence, we get Eq.\,(\ref{UE-3}) for 
the system $\Omega_\bT^W$.

Finally, we show Eq.\,(\ref{UE-4}). 
First, by Lemma \ref{key-lemma-2}, we have  
\begin{align*}
\tilde g(t)\tilde h(t)&= \lp  \sum_{T' \in \bar{\mathbb T}^W} t^{|T'|} \mathcal  V_{T'} \rp 
\lp  \sum_{T'' \in \bar{\mathbb H}^W} \beta (T'') t^{|T''|-1} 
\mathcal  V_{T''} \rp \\
&=\sum_{T \in \bar{\mathbb T}^W} t^{|T|-1} \lp
\sum_{\substack{C \in \mathcal C(T) \\ R_C(T) \in \bar{\mathbb H}^W}} 
\beta (R_C(T)) \rp \mathcal  V_T. \nno
\end{align*}

Note that, the set of all admissible 
cuts $C$ such that 
$R_C(T)  \in \bar{\mathbb H}^W$ is 
in $1$-$1$ correspondence
with the set of rooted subtrees $S$ of $T$ 
with $S\in \bar{\mathbb H}^W$. 
But any subtree 
$S$ of $T$ with $S\in \bar{\mathbb H}^W$ is 
completely determined 
by the unique leaf of $S$.
Therefore, the set of all admissible 
cuts $C$ such that 
$R_C(T)  \in \bar{\mathbb H}^W$ is 
in $1$-$1$ correspondence
with the set of non-root vertices 
of $T$. With this observation, we have, 
for any $T\in \bar{\mathbb T}^W$, 
\begin{align}
\sum_{\substack{C \in \mathcal C(T) \\ R_C(T) \in \bar{\mathbb H}^W}} 
\beta (R_C(T))=|T|.
\end{align}
Hence, combining the two equations above, we get 
$$
\tilde g(t)\tilde h(t)=\sum_{T \in \bar{\mathbb T}^W} t^{|T|-1} 
|T| \mathcal  V_T=\frac{d \tilde g(t)}{dt},
$$
which is  Eq.\,(\ref{UE-4}) for the system $\Omega_\bT^W$.
\epfv

Now, by the universal property of 
the \cNcs  system $(\cNsf$, $\Pi)$ from NCSFs
in Theorem \ref{Universal}, we have 
the following correspondence between NCSFs 
and $W$-labeled rooted trees.

\begin{theo}\label{T4.6}
For any nonempty $W\subseteq \bN^+$, 
there exists a unique homomorphism 
${\mathcal T}_W: \cNsf \to {\mathcal H}_{GL}^W$ 
of graded $K$-Hopf algebras such that 
${\mathcal T}_W^{\times 5}(\Pi)=\Omega_\bT^W$.

In particular, we have the following correspondence
from the NCSFs in $\Pi$ to the elements 
in $\Omega_\bT^W$: 
\allowdisplaybreaks{
\begin{align}
{\mathcal T}_W(\Lambda_m) &= \sum_{T \in \bar{\mathbb S}^W_m} (-1)^{o(T)+|T|}  {\mathcal  V}_T, 
\label{Lam-shrub} \\ 
{\mathcal T}_W(S_m)       &=\sum_{T \in \bar{\bT}^W_m}  {\mathcal  V_T}, \\ 
{\mathcal T}_W(\Psi_m)    &=\sum_{T \in \bar{\mathbb H}^W_m} \beta_T  \mathcal  V_T, \\ 
{\mathcal T}_W(\Phi_m)    &=m \sum_{T \in \bar{\mathbb P}_m} \theta_T {\mathcal  V_T},\\ 
{\mathcal T}_W(\Xi_m)     &=\sum_{T \in \bar{\mathbb P}_m}  \gamma_T {\mathcal  V_T}, \label{Xi-theta}
\end{align} }
for any $m\geq 1$.
\end{theo}

\pf Note that the coefficients  
of $t^m$ $(m\geq 1)$ of the generating function 
$\tilde h(t)$ (see Eq.\,(\ref{Def-3h(t)})) 
are all primitive elements of the Hopf algebra 
$\cH_{GL}^W$, since they are linear combinations 
of chains. Then, by Theorem \ref{Universal}, $(b)$, 
we have a unique homomorphism 
${\mathcal T}_W: \cNsf \to {\mathcal H}_{GL}^W$ 
of $K$-Hopf algebras such that 
${\mathcal T}_W^{\times 5}(\Pi)=\Omega_\bT^W$. 
In particular, we have 
$\cT_W (\lambda (t) )=\tilde f(t)$ 
which is same as Eq.\,(\ref{Lam-shrub}) 
for any $m\geq 1$.
Note that, both sides of Eq.\,(\ref{Lam-shrub})
have weight $m$ in $\cNsf$ and $\cH_{GL}^W$, 
respectively. 
Also note that the gradings of $\cNsf$ and $\cH_{GL}^W$ are 
given by the weights of NCSFs and $W$-labeled rooted 
trees, respectively, and $\cNsf$ is the free algebra 
generated by $\Lambda_m$ $(m\geq 1)$. By the facts above,
it is easy to check that $\cT_W$ also preserves the gradings. 
\epfv

Note that the graded duals of $\cNsf$ and $\cH_{GL}^W$ are 
the graded $K$-Hopf algebras $\cQf$ of 
quasi-symmetric functions and the Connes-Kreimer 
Hopf algebra $\cH_{CK}^W$, respectively.
Since the $K$-Hopf algebra homomorphism 
$\cT_W: \cNsf\to \cH_{GL}^W$ preserves the gradings, 
we can take the graded duals and get 
the following correspondence. 

\begin{corol}\label{cQsf}
For any non-empty $W\subseteq \bN^+$, $\mathcal T^*_W: \cH_{CK}^W\to \cQf$ is 
a homomorphism of graded $K$-Hopf algebras.  
\end{corol}

Finally, let us end this section with the following two remarks.

\begin{rmk}\label{R4.8}
As we mentioned earlier in Remark \ref{new-rmk},
by applying the specialization 
${\mathcal T}_W:\cNsf \to \cH_{GL}^W$ in Theorem \ref{T4.6},  
we will get a host of identities for the $W$-rooted 
trees on the right hand sides of 
Eqs.\,$(\ref{Lam-shrub})$--$(\ref{Xi-theta})$
from the identities of the NCSFs on the left hands.
We believe some of these identities are 
also interesting from the aspect of 
combinatorics of rooted trees. 
For example, it is not obvious at all 
that the invariants of rooted 
trees given by coefficients of 
the generating functions $\tilde f(t)$,    
$\tilde d(t)$, $\tilde h(t)$ and 
$\wtilde m(t)$ can be obtained 
by evaluating the corresponding 
$($noncommutative$)$ polynomials 
of NCSFs at coefficients of 
$\tilde g(t)$ which is just the trivial 
invariant of rooted trees 
whose value at any rooted 
tree is always $1$.

Note that the same problem 
has been studied in detail in \cite{GTS-III} 
for the differential operators
in the \cNcs over the differential operator algebras 
constructed in \cite{GTS-II}. But, in order to keep this paper 
in a certain size, we will skip the discussions on 
these identities and refer the interested reader 
to \cite{GTS-III} for similar discussions.
\end{rmk}

\begin{rmk}\label{R4.9}
In the followed paper \cite{GTS-V},
by using some relations of the \cNcs  system 
$(\cH^W_{GL}, \Omega_\bT^W)$ 
with the \cNcs  systems 
constructed in \cite{GTS-II} over 
differential operator algebras,
it will be shown that, when $W=\bN^+$, 
the $K$-Hopf algebra homomorphism 
${\mathcal T}_W:\cNsf\to \cH_{GL}^W$ 
in Theorem \ref{T4.6} is actually 
an embedding. So, in this case, 
$\cT_W^*:\cH_{CK}^W\to \cQf$ in 
Corollary \ref{cQsf} is a surjective graded 
$K$-Hopf algebra homomorphism. 

It seems that, the $K$-Hopf algebra homomorphism 
${\mathcal T}_W:\cNsf\to \cH_{GL}^W$ is injective for any non-empty
label set $W$. But we will leave this for 
future investigations.   
\end{rmk}

\renewcommand{\theequation}{\thesection.\arabic{equation}}
\renewcommand{\therema}{\thesection.\arabic{rema}}
\setcounter{equation}{0}
\setcounter{rema}{0}

\section{\bf A Combinatorial Interpretation of the Constants $\theta_T$} \label{S5}

In this section, we give a combinatorial 
interpretation for the constants $\theta_T$ 
$(T\in \bT)$ in Definition \ref{Def-theta}
which have been used in the construction 
of $\tilde d(t)$ for the \cNcs  system 
$\Omega_\bT^W$ over $\cH_{GL}^W$.
We first recall the strict order polynomials
and the order polynomials of finite posets 
(partially ordered sets) 
(see \cite{St1}).
Then, by applying some of results 
proved in \cite{WZ} for the strict 
order polynomials of rooted forests
and the well-known Reciprocity Relation 
(see Proposition \ref{Recip-Re}) between the 
strict order polynomials and the order polynomials
of finite posets, 
we show in Proposition \ref{S6-Main} 
that, for any $T\in \bar{\mathbb T}$, 
$\theta_T$ is same as the coefficient 
of $s$ of the order polynomial 
of the rooted forest $B_-(T)$. 

First, recall that, a {\it poset} is a set $P$ 
with a partial order defined for its elements.
Note that, for any unlabeled rooted forest $F\in \mathbb F$, 
the set of the vertices of $F$ has a natural partial order, 
namely, $u \leq v$ if $u=v$ or $u$ and $v$ are connected 
by a path with $v$ being further away for the root of 
the connected component. With this partial order, 
the set of the vertices of $F$ forms a finite poset, 
which we will still denote by $F$.

For any $n\geq 1$, let $[n]$ denote the totally ordered 
set $\{1, 2, ... , n\}$.  For any finite poset $P$, 
a map $f: P\to [n]$ is said to be 
{\it strictly order-preserving} (resp.\,\,{\it order-preserving}) 
if, for any $a, b\in P$ with $a<b$ in $P$, 
then $f(a)<f(b)$ (resp.\,\,$f(a)\leq f(b)$). It is well-known that, 
for each finite poset $P$, there exists a unique polynomial 
$\bar \Omega (P, s)$ (resp.\,\,$\Omega(P, s)$) in formal variable
$s$ such that, for any $n\geq 1$, 
$\bar \Omega (P, n)$ (resp.\,\,$\Omega(P, n)$) equals to the number 
of strict order-preserving (resp.\,\,{\it order-preserving}) 
from $P$ to $[n]$. The strict order polynomial 
and the order polynomial of finite posets 
are related by the so-called reciprocity relation.

\begin{propo} \label{Recip-Re}
$(\text{\bf Reciprocity Relation})$
For any fixed finite poset, we have
\begin{align}\label{Recip-Re-e1}
 \Omega(P, s)=(-1)^{|P|}\bar \Omega(P, -s),
\end{align}
where $|P|$ denotes the cardinal number of the finite set $P$.
\end{propo}

For a proof of this remarkable result, 
see Corollary $4.5.15$ in \cite{St1}. 

By Corollary $3.15$ in \cite{WZ} (also see Example $4.4$ there) 
and Eq.\,(\ref{Recip-Re-e1}) above, we get the (strict) order 
polynomials of chains and shrubs as follows.

\begin{exam}\label{Exam6.2} 
For any $m\geq 1$, let $C_m$ be the chain of hight $m-1$
and $S_m$ the shrub with $m$ leaves. Then we have
\begin{align}
\bar \Omega (C_m, s)&=\binom{s}{m}=\frac {s(s-1)\cdots (s-m+1)}{m!},\label{Exam6.2-e1} \\
\Omega (C_m, s)&=(-1)^m \binom{-s}{m}=\binom{s+m-1}{m},\label{Exam6.2-e2} \\
\bar \Omega (S_m, s)&=\int_0^s B_m (u)\,du,
=\frac{B_{m+1}(s)-B_{m+1}(0)}{m+1}, \label{Exam6.2-e3} \\
\Omega (S_m, s)&=(-1)^{m+1} \int_0^{-s} B_m (u)\,du, \label{Exam6.2-e4}
\end{align}
where $B_m(u)$ $(m\geq 1)$ are the
Bernoulli polynomials which are defined by the generating function
$$
\frac{xe^{ux}}{e^{x}-1}=1 + \sum_{m=1}^{\infty}B_m(u)\frac{x^{m}}{m!}.
$$
\end{exam}

\begin{propo}\label{P6.3}
$(a)$ For any finite poset $P$, $\bar \Omega (P, 0)=0$.


$(b)$ For any rooted forest 
$F=T_1T_2\cdots T_m$ with $T_i\in \mathbb T$ $(1 \leq i \leq m)$, 
we have
\begin{align}\label{P6.3-e1}
\bar \Omega(F, s)=\bar \Omega(T_1, s)\bar \Omega(T_2, s) \cdots \bar \Omega(T_m, s).
\end{align}

$(c)$ 
For any unlabeled rooted tree $T$ 
with $T=B_+(F)$, we have
\begin{align}
\Delta \bar \Omega(T, s)&=\bar \Omega(F, s),\label{P6.3-e2}\\
\nabla \Omega(T, s)&=\Omega(F, s),\label{P6.3-e3}
\end{align}
where $\Delta$ and $\nabla$ are the linear operators from 
$K[s]\to K[s]$ defined by, for any $f(s)\in K[s]$,
$\Delta f(s)=f(s+1)-f(s)$ and 
$\nabla f(s)=f(s)-f(s-1)$, 
respectively.
\end{propo}
 
$(a)$ is well-known, for example, it can be easily proved
by the recurrent formulas in \cite{SWZ}
for the (strict) order polynomials. $(b)$
follows directly from the definition of 
the strict order polynomials. $(c)$ was first proved 
by J. Shareshian (unpublished). 
For a proof of Eq.\,(\ref{P6.3-e2}), 
see Theorem $4.5$ in \cite{WZ}. 
Eq.\,(\ref{P6.3-e3}) can be proved 
similarly as Eq.\,(\ref{P6.3-e2}).
For more studies on these properties 
of the (strict) order polynomials, 
see \cite{Tree-Inv} and \cite{SWZ}.

Now, for any finite poset $P$, we define 
\begin{align}
\phi_P &:=\left. \frac{d}{d s} \bar \Omega(P, s)\right |_{s=0},\\
\varphi_P& :=\left. \frac{d}{d s}  \Omega(P, s) \right |_{s=0}.\label{Def-varphi-P}
\end{align}

By Proposition \ref{Recip-Re} and \ref{P6.3}, it is easy to see that 
we have the following corollary.

\begin{corol}\label{C6.4}
$(a)$ For any finite poset $P$, we have $\phi_P=(-1)^{|P|-1}\varphi_P$.
 
$(b)$ For any rooted forest $F\in \mathbb F$ with more than one connected component, 
we have $\phi_F= \varphi_F=0$.
\end{corol}

The following proposition have been proved in \cite{WZ}. But note that 
the definition of $T_{\vec e,k}$ we adapt here is different 
form the one used in \cite{WZ}. So the equations 
in the proposition below
have been modified accordingly.

\begin{propo}\label{P6.5}
$(a)$ The constants $\{\phi_T\,|\, T \in \mathbb T \}$
 satisfy, and are uniquely determined by
\begin{align*}
& \phi_{T=\circ} =1, \\
& \phi_T =-\sum_{k=2}^{v(T)} \frac 1{k!}
\sum_{\substack{\vec e=(e_1, \ldots, e_{k-1})\in E(T)^{k-1}\\e_1\succ
\cdots \succ e_{k-1}}}
\phi_{B_-(T_{\vec e,1})}  \cdots  \phi_{B_-(T_{\vec e, k-1})} \phi_{T_{\vec e,k}}, \nno
\end{align*}
when $v(T) \geq 2$. 

$(b)$ For any $T\in \mathbb T$, we have
\begin{align*}
\bar \Omega(T, s)=\phi_T \, s +\sum_{k=2}^{v(T)} \frac {s^{k}}{k!}
\sum_{\substack{\vec e=(e_1,\ldots,e_{k-1})\in E(T)^{k-1}\\e_1\succ
\cdots\succ e_{k-1}}}
\phi_{B_-(T_{\vec e,1})}  \cdots  \phi_{B_-(T_{\vec e, k-1})} \phi_{T_{\vec e,k}}.
\end{align*}
\end{propo}

By Proposition \ref{P6.5}, 
Corollary \ref{C6.4}, $(a)$ and Proposition \ref{Recip-Re}, 
it is easy to check that 
the order polynomials $\Omega(T, s)$ of rooted trees $T$
can be obtained as follows. 

\begin{corol}\label{Recur-varphi}
$(a)$ The constants $\{\varphi_T | T \in \mathbb T \}$
 satisfy and are uniquely determined by
\begin{align*}
& \varphi_{T=\circ} =1,\\
& \varphi_T = \sum_{k=2}^{v(T)} \frac {(-1)^k}{k!}
\sum_{\substack{\vec e=(e_1, \ldots, e_{k-1})\in E(T)^{k-1}\\e_1\succ
\cdots \succ e_{k-1}}}
\varphi_{B_-(T_{\vec e,1})} \cdots  \varphi_{B_-(T_{\vec e,k-1})} \varphi_{T_{\vec e,k}} \nno
\end{align*}
when $v(T) \geq 2$. 

$(b)$ For any $T\in \mathbb T$, we have
\begin{align}\label{Recur-varphi-e2}
\Omega(T, s)= \varphi_T\,  s + \sum_{k=2}^{v(T)} \frac {s^{k}}{k!}
\sum_{\substack{\vec e=(e_1,\ldots,e_{k-1})\in E(T)^{k-1}\\e_1\succ
\cdots\succ e_{k-1}}}
\varphi_{B_-(T_{\vec e,1})} \cdots  \varphi_{B_-(T_{\vec e,k-1})} \varphi_{T_{\vec e,k}}.
\end{align}
\end{corol}

Note that, from the definition of the order polynomials, 
we have $\Omega(P, 1)=1$ for any finite poset $P$. Using this fact
and evaluating $\Omega(T, s)$ in Eq.\,(\ref{Recur-varphi-e2}) 
at $s=1$, 
we get another recurrent formula 
for the constants $\varphi_T$ 
of rooted trees.

\begin{corol} \label{2nd-Recur-varphi}
 The constants $\{\varphi_T | T \in \mathbb T \}$
 satisfy and are uniquely determined by
\begin{align}
& \varphi_{T=\circ} =1, \nno \\
& \varphi_T = 1- \sum_{k=2}^{v(T)} \frac {1}{k!}
\sum_{\substack{\vec e=(e_1,\ldots,e_{k-1})\in E(T)^{k-1}\\e_1\succ
\cdots\succ e_{k-1}}}
\varphi_{B_-(T_{\vec e,1})} \cdots  \varphi_{B_-(T_{\vec e,k-1})} 
\varphi_{T_{\vec e,k}}. \label{2nd-Recur-varphi-e1}
\end{align}
\end{corol}

Now, we consider the constants $\theta_T$ $(T\in \mathbb T)$
defined in Definition \ref{Def-theta} 
in Section \ref{S4} and prove the following main result
of this section.  

\begin{propo}\label{S6-Main}
For any $T\in {\mathbb T}$ with $T=B_+(F)$,  we have
\begin{align}
\theta_T = \varphi_F. \label{S6-Main-e1}
\end{align}
\end{propo}
\pf For any $T\in {\mathbb T}$ with $T=B_+(F)$, 
we set $\tilde \theta_T:= \varphi_F$. For convenience, 
we also set $\tilde \theta_\circ=0$. We need show that
$\tilde \theta_T=\theta_T$ for any $T\in {\mathbb T}$.

Note first that, if $T$ the singleton $\circ $ or a 
non-primitive rooted tree, i.e. $F$ is the empty or 
has at least two connected components,
by the definition of $\theta_T$ in 
Definition \ref{Def-theta} and 
Corollary \ref{C6.4}, $(b)$,
we have $\theta_T=\tilde \theta_T=0$.
To show that $\tilde \theta_T=\theta_T$ 
for all primitive rooted trees $T$, it will 
be enough to show that, $\tilde \theta_T$ 
$(T\in \mathbb P)$
also satisfies the recurrent relations 
in Definition \ref{Def-theta}.

We use the mathematical induction on $v(T)$. 
First, for the case $v(T)$, i.e.
$T=B_+(\circ)=C_2$.
Since $\Omega(\circ, s)=s$,  
$\tilde \theta_T=\varphi_\circ = 1$.
While $\theta_T$ is defined to be $1$ in 
Definition \ref{Def-theta}. Hence, 
 $\theta_T=\tilde \theta_T$ in this case.

Now, assume  
$T=B_+(T')$ with $T'\in \mathbb T$ and $v(T')\geq 2$. 
Applying Eq.\,(\ref{2nd-Recur-varphi-e1}) to $T'$, we have
\begin{align*}
& \varphi_{T'}= 1 - \sum_{k=2}^{v(T')} \frac 1{k!}
\sum_{\substack{\vec e=(e_1, \ldots, e_{k-1})\in E(T')^{k-1}\\e_1\succ
\cdots \succ e_{k-1}}}
\varphi_{B_-(T_{\vec e,1})}  \cdots 
\varphi_{B_-(T_{\vec e, k-1})} \varphi_{T_{\vec e,k}} \, . \nno
\end{align*}

Note that, we can identify the set of
$\vec e=(e_1, \ldots, e_{k-1})\in E(T')^{k-1}$  
with the set of $\vec e=(e_1, \ldots, e_{k-1})\in E(T)^{k-1}$ 
such that $T_{\vec e,k}\neq \circ$. 
With these observations, by replacing the constants $\varphi$'s 
by $\tilde \theta$'s
in the summation of the equation above, we have 

\begin{align*}
\tilde \theta_T &= 1 - \sum_{k=2}^{v(T)} \frac 1{k!}
\sum_{\substack{\vec e=(e_1, \ldots, e_{k-1})\in E(T)^{k-1}\\e_1\succ
\cdots \succ e_{k-1} \\
T_{\vec e,k}\neq \circ }}
\tilde \theta_{T_{\vec e,1}} \tilde \theta_{T_{\vec e,2}} \cdots  \tilde \theta_{T_{\vec e,k}}\\
\intertext{Applying the induction assumption to $T_{\vec e,j}$'s and using
the fact that $\tilde \theta_\circ = 0$:}
&=1 - \sum_{k=2}^{v(T)} \frac 1{k!}
\sum_{\substack{\vec e=(e_1, \ldots, e_{k-1})\in E(T)^{k-1}\\e_1\succ
\cdots \succ e_{k-1} }}
\tilde \theta_{T_{\vec e,1}} \tilde \theta_{T_{\vec e,2}} 
\cdots  \tilde \theta_{T_{\vec e,k}}.\\
\end{align*}

Therefore, the constants $\{ \tilde\theta_T \,|\, T\in \mathbb P \}$ 
also satisfy the recurrent relations of 
$\{\theta_T \,|\, T\in \mathbb P \}$ 
in Definition \ref{Def-theta}.
Hence, we have $\theta_T=\tilde \theta_T=\varphi_T$ 
for any $T\in \mathbb T$. 
\epfv

For an interpretation of the constant
$\phi_T$, which we have shown is same as 
$(-1)^{v(T)-1}\varphi_T=(-1)^{v(T)-1}\theta_{B_+(T)}$,  
in terms of the numbers of chains with fixed lengths in  
the lattice of the ideals of the poset $T$, 
see Lemma $2.8$ in \cite{SWZ}.

\begin{corol}\label{C6.10}
For any $T\in \mathbb P$, we have
\begin{align}\label{C6.10-e1}
\nabla \Omega(T, s) & = \theta_T \,  s +  \sum_{k=2}^{v(T)} \frac {s^{k}}{k!}
\sum_{\substack{\vec e=(e_1,\ldots,e_{k-1})\in E(T)^{k-1}\\e_1\succ
\cdots\succ e_{k-1}}}
\theta_{B_-(T_{\vec e,1})}  \cdots \theta_{B_-(T_{\vec e,k-1})} \theta_{T_{\vec e,k}}.
\end{align}
In particular, $\theta_T$ is also  
the coefficient of $s$ of the polynomial $\nabla \Omega(T, s)$.
\end{corol}

\pf  First, we write $T=B_+(T')$ with $T'\in \bT$. By Eq.\,$(\ref{P6.3-e3})$ 
and Proposition \ref{P6.5}, $(b)$, we have
\begin{align*}
\nabla &\Omega(T, s)=\Omega(T', s)\\
&=\varphi_{T'} \,  s + \sum_{k=2}^{v(T')} \frac {s^{k}}{k!}
\sum_{\substack{\vec e=(e_1,\ldots,e_{k-1})\in E(T')^{k-1}\\e_1\succ
\cdots\succ e_{k-1}}}
\varphi_{B_-(T_{\vec e,1})} \cdots \varphi_{B_-(T_{\vec e,k-1})} 
\varphi_{T'_{\vec e,k}}.
\end{align*}

Then, applying Eq.\,(\ref{S6-Main-e1}) and replacing 
the constants $\varphi$'s in the sum above by the constant
$\theta$'s, we get Eq.\,(\ref{C6.10-e1}).
\epfv

{\small \sc Department of Mathematics, Illinois State University,
Normal, IL 61790-4520.}

{\em E-mail}: wzhao@ilstu.edu.

\end{document}